\input form
\font\teneufm=eufm10
\font\seveneufm=eufm7
\font\fiveeufm=eufm5
\newfam\eufmfam
\textfont\eufmfam=\teneufm
\scriptfont\eufmfam=\seveneufm
\scriptscriptfont\eufmfam=\fiveeufm


\def\uno{{\mat 1}}
\def\Rmatw{\Rmat_\omega}
\def\Dmatw{\Dmat_\omega}
\def\Lambdamatw{\Lambdamat_\omega}

\def\frac#1#2{{#1 \over #2}}
\def\parder#1#2{{\partial#1 \over \partial#2}}

\def\build#1_#2^#3{\mathrel{
\mathop{\kern 0pt#1}\limits_{#2}^{#3}}}

\def\lie#1{L_{#1}}
\def\liec#1{L_{#1}}

\def\ordnorma#1{\|{#1}\|}
\def\bignorma#1{\bigl\|{#1}\bigr\|}

\def\biggnorma#1{\biggl\|{#1}\biggr\|}

\def\projn{P_{{\cal N}}}
\def\projr{P_{{\cal R}}}
\def\chr{{\it X$\rho$\'o$\nu o\varsigma$}}

\newdimen\pagewidth \newdimen\pageheight \newdimen\ruleht
\pagewidth=\hsize \pageheight=\vsize \ruleht=.5pt
\newdimen\halfwd \halfwd=0.48\pagewidth



\cita{Bazzani-1988}{A. Bazzani, P. Mazzanti, G. Servizi, G.
Turchetti: {\it Normal forms for Hamiltonian maps and nonlinear
effects in a particle accelerator}, Nuovo Cim. B, {\bf 102}, 51--80
(1988).}

\cita{Bazzani-1993}{A. Bazzani, M. Giovannozzi, G. Servizi, E.
Todesco, G. Turchetti: {\it Resonant normal forms, interpolating
Hamiltonians and stability of area preserving maps}, Physica D, {\bf
64}, 66--93 (1993).}

\cita{Bazzani-1994}{A. Bazzani, E. Todesco, G. Turchetti: {\it A
normal form approach to the theory of nonlinear betatronic motion},
CERN Report No.94-02 (1994).}

\cita{BenGalGio-1985}{G. Benettin, L. Galgani, A. Giorgilli: {\it A proof of 
Nekhoroshev's theorem for the stability times in nearly integrable
Hamiltonian systems.} Cel. Mech., {\bf 37}, 1--25 (1985).}

\cita{Birkhoff-1920}{G.D. Birkhoff: {\it Surface transformations and
thair dynamical applications}, Acta Mathematica, {\bf 43}, 1--119
(1920).}

\cita{Birkhoff-1927}{G.D. Birkhoff: {\it Dynamical systems}, New York (1927).}

\cita{BorCarSkoPapVit-2012}{J. Boreux, T. Carletti, C. Skokos, Y. Papaphilippou,
M. Vittot: {\it Efficient control of accelerator maps}, International
Journal of Bifurcation and Chaos, {\bf 22} (9), 1250219-1--1250219-9
(2012).}

\cita{BorCarSkoVit-2012}{J. Boreux, T. Carletti, C. Skokos, M. Vittot: 
{\it Hamiltonian control used to improve the beam stability in
particle accelerator models}, Communications in Nonlinear Science and
Numerical Simulation, {\bf 17}, 1725--1738 (2012).}

\cita{BouSko-2006}{T. Bountis and C. Skokos: {\it Space charges can
significantly affect the dynamics of accelerator maps}, Phys. Lett. A,
{\bf 358}, 126--133 (2006).}

\cita{CalDes-1991}{F. Callier and C. Desoer: {\it Linear system
theory}, Springer-Verlag, (1991).}

\cita{Celletti-1988}{A. Celletti and L. Chierchia: {\it
Construction of analytic KAM surfaces and effective stability bounds},
Communications in Mathematical Physics, {\bf 118}, 119--161 (1988).}

\cita{Celletti-1997}{A. Celletti, L. Chierchia: {\it On the
stability of realistic three-body problems}, Communications in
Mathematical Physics, {\bf 186}, 413--449 (1997).}


\cita{ChaVitElsCirPet-2005}{C. Chandre, M. Vittot, Y. Elskens,
G. Ciraolo, M. Pettini: {\it Controlling
chaos in area-preserving maps}, Physica D {\bf 208}, 131--146 (2005).}

\cita{ChaCirDovLimMacVit-2005}{C. Chandre, G. Ciraolo, F. Doveil, R. Lima, A. Macor, M. Vittot:
{\it Channeling chaos by building barriers}, Phys. Rev. Letters {\bf
 94}, 074101 (2005).}

\cita{Cherry-1924a}{T.M. Cherry: {\it On integrals developable about a 
singular point of a Hamiltonian system of differential equations},
Proc.  Camb. Phil. Soc., {\bf 22}, 325--349 (1924).}

\cita{Cherry-1924b}{T.M. Cherry: {\it On integrals developable about a 
singular point of a Hamiltonian system of differential equations, II}, Proc. 
Camb. Phil. Soc., {\bf 22} 510--533 (1924).
}

\cita{CirBriChaFloLimVitPet-2004}{G. Ciraolo, F. Briolle, C. Chandre, E. Floriani,
R. Lima, M. Vittot, M. Pettini: {\it Control of Hamiltonian chaos as a
possible tool to control anomalous transport in fusion plasmas}
Physical Review E, {\bf 69} (5), 056213 (2004).}

\cita{CirChaLimVitPet-2004}{G. Ciraolo, C. Chandre, R. Lima, M. Vittot, 
M. Pettini: {\it Control of chaos in Hamiltonian systems} CeMDA, {\bf
90}, 3--12 (2004).}

\cita{CirChaLimVitPetFigGhe-2004}{G. Ciraolo, C. Chandre, R. Lima, M. Vittot, M. Pettini,
C. Figarella, P. Ghendrih: {\it Controlling chaotic transport in a
Hamiltonian model of interest to magnetized plasmas}, Journal of
Physics A {\bf 37}, 3589 (2004).}

\cita{Deprit-1969}{A. Deprit: {\it Canonical transformations
depending on a small parameter}, {Cel. Mech.}, {\bf 1}, 12--30 (1969).}

\cita{Fasso-1989}{F. Fass\`o: {\it On a relation among Lie series
}, Cel. Mech., {\bf 46}, 113--118 (1989).}

\cita{Gallavotti-1982}{G. Gallavotti: {\it A criterion of of
integrability for perturbed nonresonant harmonic oscillators. ``Wick
ordering'' of the perturbations in Calssical Mechanics and invariance
in the frequency spectrum}, Comm. Math. Phys., {\bf 87}, 365--383
(1982).}

\cita{GelGel-2009}{V. Gelfreich and N. Gelfreikh: {\it Unique resonant
normal form for area preserving maps at an elliptic fixed point},
Nonlinearity {\bf 22}, 783--810 (2009).}

\cita{GioGal-1978}{A. Giorgilli and L. Galgani: {\it Formal
integrals for an autonomous Hamiltonian system near an equilibrium
point}, Cel. Mech., {\bf 17}, 267--280 (1978).}

\cita{Giorgilli-1985}{A. Giorgilli and L. Galgani: {\it Rigorous
estimates for the series expansions of Hamiltonian perturbation
theory}, Cel. Mech., {\bf 37}, 95--112 (1985).}

\cita{GioDelFinGalSim-1989}{A. Giorgilli, A. Delshams, E. Fontich, L. Galgani, C. 
 Sim\'o: {\it Effective stability for a Hamiltonian system near an
 elliptic equilibrium point, with an application to the restricted
 three body problem.} J. Diff. Eqs., {\bf 77}, 167--198 (1989).  }

\cita{GioZeh-1992}{A. Giorgilli and E. Zehnder: {\it 
Exponential stability for time dependent potentials}, ZAMP, {\bf 5},
827--855 (1992).}

\cita{GioSko-1997}{A. Giorgilli and Ch. Skokos: {\it On the
stability of the Trojan asteroids}, Astron. Astroph., {\bf 317},
254--261 (1997).}

\cita{Giorgilli-2003}{A. Giorgilli: {\it Notes on exponential
    stability of Hamiltonian systems}, in {\it Dynamical Systems, Part
    I\/}. Pubbl. Cent. Ric. Mat. Ennio De Giorgi, Sc. Norm. Sup.
    Pisa, 87--198 (2003).}

\cita{GioLocSan-2009}{A. Giorgilli, U. Locatelli, M. Sansottera:
{\it Kolmogorov and Nekhoroshev theory for the problem of three
bodies}, CeMDA, {\bf 104}, 159--175 (2009).}

\cita{Giorgilli-2013}{A. Giorgilli: {\it On the representation of
maps by Lie transforms}, Rendiconti dell'Istituto Lombardo Accademia
di Scienze e Lettere, Classe di Scienze Ma\-te\-ma\-ti\-che e Naturali, {\bf
146}, 251--277 (2012).}


\cita{GioSan-2012}{A. Giorgilli and  M. Sansottera: {\it Methods of algebraic manipulation in
    perturbation theory}, Workshop Series of the Asociacion Argentina
  de Astronomia, {\bf 3}, 147--183 (2011).}

\cita{Groebner-1957}{W.\ Gr\"obner: {\it Nuovi contributi alla 
teoria dei sistemi di equazioni differenziali nel campo analitico}, 
Atti Accad. Naz. Lincei. Rend. Cl. Sci. Fis. Mat. Nat., {\bf 23},
375--379 (1957).}

\cita{Groebner-1967}{W. Gr\"obner: {\it Die  Lie-Reihen und Ihre 
Anwendungen}, VEB Deutscher Verlag der Wissenschaften, Berlin (1967).}

\cita{HanCar-1984}{J. Hanson and J. Cary: {\it Elimination of stochasticity in stellerators}, Physics of fluids, {\bf 27}, 767--769, (1984).}

\cita{HenRoe-1974}{J. Henrard and J. Roels: {\it Equivalence for Lie
transforms}, Cel. Mech., {\bf 10}, 497--512 (1974).}

\cita{HinPri-2005}{D. Hinrichsen and A. Pritchard: {\it Mathematical systems theory I modelling, state space analysis, stability and robustness}, Springer-Verlag, (2005).}

\cita{Hori-1966}{G. Hori: {\it Theory of general perturbations with
unspecified canonical variables}, Publ. Astron. Soc. Japan, {\bf 18},
287--296 (1966).}

\cita{Kolmogorov-1954}{A.N. Kolmogorov: {\it Preservation of
conditionally periodic movements with small change in the Hamilton
function}, Dokl. Akad. Nauk SSSR, {\bf 98}, 527 (1954). English
translation in: Los Alamos Scientific Laboratory translation
LA-TR-71-67; reprinted in: G.\ Casati, J.\ Ford: Stochastic behavior
in classical and quantum Hamiltonian systems, Lecture Notes in Physics,
{\bf 93}, 51--56 (1979).}

\cita{Littlewood-1959a}{J.E. Littlewood: {\it On the equilateral 
configuration in the restricted problem of three bodies}, Proc. London 
Math. Soc.(3), {\bf 9}, 343--372 (1959).}

\cita{Littlewood-1959b}{J.E. Littlewood: {\it The Lagrange configuration 
in celestial mechanics}, Proc. London Math. Soc.(3), {\bf 9},
525--543 (1959).}

\cita{Locatelli-2001}{A. Locatelli: {\it Optimal control: an introduction.}, Birkha\"user (2001).}

\cita{LocGio-2000}{U. Locatelli and A. Giorgilli: {\it  Invariant
tori in the secular motions of the three-body planetary systems},
Cel. Mech., {\bf 78}, 47--74 (2000).}  

\cita{LocGio-2007}{U. Locatelli and A. Giorgilli: {\it Invariant
tori in the Sun-Jupiter-Saturn system}, DCDS-B, {\bf 7}, 377--398 (2007).} 

\cita{MorGio-1995a}{A. Morbidelli and A. Giorgilli: {\it Superexponential 
stability of KAM tori}, J. Stat. Phys., {\bf 78}, 1607--1617 (1995).}

\cita{MorGio-1995b}{A. Morbidelli and A. Giorgilli: {\it On a connection 
between KAM and Nekhoroshev's theorem}, Physica D, {\bf 86}, 514--516 
(1995).}

\cita{Moser-1955}{J. Moser: {\it Stabilit\"atsverhalten kanonisher 
differentialgleichungssysteme,} Nachr.\ Akad.\ Wiss.\ G\"ottingen,
Math.\ Phys.\ K1.\ IIa, {\bf 6} 87--120 (1955).}


\cita{Nekhoroshev-1977}{N.N. Nekhoroshev: {\it Exponential estimates of the stability
time of near-integrable Hamiltonian systems}. English translation:
Russ. Math.  Surveys, {\bf 32}, 1 (1977).  }
 
\cita{Nekhoroshev-1979}{N.N. Nekhoroshev: {\it Exponential estimates of the stability
time of near-integrable Hamiltonian systems, 2.} Trudy Sem.\ Im.\ G.\
Petrovskogo, {\bf 5}, 5 (1979).  English translation: {\it Topics in
modern Mathematics}, Petrovskij Semin., {\bf 5}, 1--58 (1985).}

\cita{Poincare-1890}{H. Poincar\'e: {\it Sur le probl\`eme des trois
corps et les \'equations de la dynamique},  Acta Mathematica (1890).}

\cita{Poincare-1892}{H. Poincar\'e: {\it Les m\'ethodes nouvelles de la
m\'ecanique c\'eleste}, Gau\-thier-Villars, Paris (1892).}

\cita{Poincare-1912}{H. Poincar\'e: {\it Sur un th\'eor\`eme en
g\'eom\'etrie}, Rendiconti del Circolo Matematico di Palermo, {\bf
33}, 375--40 (1912).}

\cita{Raubenheimer-1993}{T. Raubenheimer: {\it The preservation of low
emittance flat beams}, SLAC-PUB-6117, 1--5 (1993).}

\cita{SanLocGio-2010}{M. Sansottera, U. Locatelli and A. Giorgilli: 
{\it On the stability of the secular evolution of the planar
Sun-Jupiter-Saturn-Uranus system}, Math. Comput. Simul. {\bf 88},
1--14 (2013).}

\cita{ScaTur-1990}{W. Scandale and G. Turchetti (Editors):
{\it Nonlinear problems in future particle accelerators}, World
Scientific, Singapore (1990).}

\cita{Takens-1974}{F. Takens: {\it Forced oscillations and bifurcations.
Applications of global analysis, I}
Commun. Math. Inst. Rijksuniv. Utrecht, {\bf 3}, (1974); reprinted in:
H.W. Broer, B. Krauskopf, G. Vegter: {\it Global
Analysis of Dynamical System}, IoP (2001).}

\cita{Todesco-1994}{E. Todesco: {\it Analysis of resonant structures
of four-dimensional symplectic mappings, using normal forms}, Physical
Review E, {\bf 50}, 4298--4301 (1994).}

\cita{Vittot-2004}{M. Vittot: {Perturbation theory and control in
classical or quantum mechanics by an inversion formula}, J. of Physics
A {\bf 37}, 6337--6357 (2004).}

\cita{VitChaCirLim-2005}{M. Vittot, C. Chandre, G. Ciraolo, R. Lima:
{\it Localised control for non-resonant Hamiltonian systems},
 Nonlinearity {\bf 18} 423--440, (2005).}

\cita{Whittaker-1916}{E.T. Whittaker: {\it On the adelphic integral of the 
differential equations of dynamics}, Proc. Roy Soc. Edinburgh, Sect. A, 
{\bf 37}, 95--109 (1916).}

\cita{WanCarSha-1996}{W. Wan, J. R. Cary, S. G. Shasharina: {\it A method for finding 4D symplectic maps with reduced chaos}, Epac 96: Proceedings of the Fifth European Particle Accelerator Conference, Sitges (Barcelona), 10--14 June 1996.}

\cita{WanCar-2001}{W. Wan and J. R. Cary:{\it Method for enlarging the dynamic aperture of accelerator lattices}, Physical Review Special Topics - Accelerators and Beams, {\bf 4}, 084001-1--084001-10 (2001).}


\def\testos{\hfill M. Sansottera, A. Giorgilli and T. Carletti}
\def\testod{High-order control for symplectic maps}

\title{High-order control for symplectic maps}
\author{\it MARCO SANSOTTERA\hfil\break 
Dipartimento di Matematica, Universit\`a degli Studi di
Milano,\hfil\break via Saldini 50, 20133\ ---\ Milano, Italy.}

\author{\it ANTONIO GIORGILLI\hfil\break
Dipartimento di Matematica, Universit\`a degli Studi di Milano,\hfil\break 
via Saldini 50, 20133\ ---\ Milano, Italy.}

\author{\it TIMOTEO CARLETTI\hfil\break
Departement of Mathematics and\hfil\break
Namur Center for Complex Systems - naXys, University of Namur,\hfil\break
8 Rempart de la Vierge, B5000\ ---\ Namur, Belgium.}

\abstract{
We revisit the problem of introducing an {\corsivo a priori control}
for devices that can be modeled via a symplectic map in a neighborhood
of an elliptic equilibrium.  Using a technique based on Lie transform
methods we produce a normal form algorithm that avoids the usual step
of interpolating the map with a flow.  The formal algorithm is
completed with quantitative estimates that bring into evidence the
asymptotic character of the normal form transformation.  Then we
perform an heuristic analysis of the dynamical behavior of the map
using the invariant function for the normalized map.  Finally, we
discuss how control terms of different orders may be introduced so as
to increase the size of the stable domain of the map.  The numerical
examples are worked out on a two dimensional map of H\'enon type.}

\section{sec:intro}{Introduction}
To control engineered human devices is a necessity to ensure their
right behavior, that is the system should follow as close as possible
the wanted trajectory, regardless from any deviations due to noise
and possible errors. For this reason the so called {\corsivo problem of
control} of dynamical systems is a long lasting research field where
engineers, physicists and mathematicians, among others, have been very
active.

Nowadays this represents a whole field in applied mathematics and in
engineering with sub-fields such as control and
systems\bibref{HinPri-2005}, linear
systems\bibref{CalDes-1991} or optimal
control\bibref{Locatelli-2001}. Because of the advancement of
technology, controllers are in action almost everywhere in our daily routine,
from the design of the water tank of the ordinary flush toilet to the
design of a lateral and longitudinal control of a Boeing (and thus the
autopilot) or in the satellite’s attitude control.

The common feature of most of the above quoted theories lies in the
idea of {\corsivo feedback}, that is one should access to the actual
state of the system (or to any relevant measure of it) and then act on
the system to achieve the desired goal, usually to stabilize the
trajectory around the nominal one. The controller is thus switched on
and off when required {\corsivo during} the evolution of the
system. This implicitly means that the typical time scale of the
system dynamics is several order of magnitude bigger than the time
required to the controller to read the state of the system, to compute
the correction and eventually modify the system accordingly.

A procedure based on feedback is clearly a hard task, if not
impossible, when one deals with very complex devices such as,
e.g., a particle accelerator, mainly because the previous assumption
about the time scale is no longer satisfied: no control could react as
fast as the particles moving at almost the speed of light.  In this
case however a control can be obtained by some external tuning to be
performed {\corsivo before} the device is turned on.  We hereby develop
this idea by presenting a method that can be applied in case
a reliable mathematical model of the device is available.  The basis
of the method is to introduce high order controllers without the need
of switches and based on the identification of {\corsivo dangerous
terms}, that once removed - or they impact reduced - will allow to
achieve the desired goal; the method is thus, in some sense, based on
an {\corsivo a priori} analysis.

The above idea has been exploited already in the end of the eighties:
see the volume~\dbiref{ScaTur-1990}.  Recently the method has been
widely investigated on the basis of the ideas presented by
Vittot\bibref{Vittot-2004} in the continuous time case and then
adapted to the discrete time case, namely maps (see,
e.g.,~\dbiref{ChaVitElsCirPet-2005}, \dbiref{VitChaCirLim-2005}
and \dbiref{BorCarSkoVit-2012}).  The introduction and the intensive
use of the Lie transform methods is the key ingredient to easily
obtain high order controllers.

For a sake of clarity we consider the problem of control in the
Hamiltonian framework, let us stress however that the following ideas
can be straightforwardly extended beyond it. In its simplest
formulation it can thus be stated as follows.  Consider a nearly
integrable canonical system of differential equations in the
neighborhood of an elliptic equilibrium, as described by the
Hamiltonian
$$
H(x,y) = H_0(x,y) + F(x,y)\ ,\quad
H_0(x,y) = \frac{1}{2} \sum_{l=1}^{n} \omega_l(x_l^2+y_l^2)
\formula{eq:1}
$$
where $(x,y)\in\reali^{2n}$ are the canonically conjugate coordinates,
$\omega\in\reali^n$ are the frequencies and $F(x,y)$ is either a power
series or a polynomial of finite degree starting with terms of degree
at least 3.  The problem is to find a function $G(x,y)$ such that the
Hamiltonian $H_0+F-G$ can be put in a simpler form, in most
applications this would mean to conjugate it to an integrable one by a
canonical transformation.  The answer should be non-trivial, in the
sense that the obvious choice $G=F$ is not accepted, because the
nonlinear character of the system should be preserved in some form.
E.g., one looks for a nonlinear integrable Hamiltonian.  It is rather
requested that $G$ should be smaller than $F$, e.g., in some norm.

A very similar problem is concerned with symplectic maps in a
neighborhood of an elliptic equilibrium.  One considers a map of a
neighborhood of the origin in the plane $\reali^2$ that is written as
$$
\pmatrix{x' \cr y'\cr} = 
 \Lambdamatw \pmatrix{x\cr y\cr} + \pmatrix{f(x,y)\cr g(x,y)\cr} 
\formula{map:1}
$$
where $\Lambdamatw$ is a suitable rotation matrix (see
formula~\frmref{eq:lambda}), while $f(x,y)$ and $g(x,y)$ are either
power series or polynomials starting with terms at least of degree 2,
and are required to satisfy the symplecticity conditions.  Again the
problem is to add a non-trivial control term such that
the resulting modified map is conjugated to a rotation, possibly a
twist one.

In the present paper we revisit the problem of control referring in
particular to the case of maps.  We reformulate it using the tool of
Lie transforms, and point out different ways of introducing control
terms.

It should be noticed that the problem of control is essentially a
particular formulation of the classical ``general problem of
dynamics'', so named by Poincar\'e (see~\dbiref{Poincare-1892},
Vol.~I,~\S$\,$13).  That is, to investigate the dynamics of a
Hamiltonian system
$$
H(p,q) = h(p) + \epsilon f(p,q;\epsilon)
\formula{eq:2}
$$
where $(p,q)\in\Gscr\times\toro^n$ with $\Gscr\subset\reali^n$ are
action-angle variables and $\epsilon$ is a small parameter.  The
Hamiltonian is assumed to be holomorphic in all variables and in the
parameter.  As proved by Poincar\'e himself, the system is generically
non-integrable due to the extreme complexity of the orbits.  However
some insight on the dynamics of the system is available by searching
for weaker properties than complete integrability.  Just to quote some
results that came after Poincar\'e, we mention the existence of
invariant tori\bibref{Kolmogorov-1954}, the theory of exponential
stability\bibref{Littlewood-1959a}\bibref{Littlewood-1959b}\bibref{Moser-1955}\bibref{Nekhoroshev-1977}\bibref{Nekhoroshev-1979}\bibref{BenGalGio-1985}, the super-exponential stability\bibref{MorGio-1995a}\bibref{MorGio-1995b}.  After
Poincar\'e's work it was soon remarked that the problem is
substantially simplified if one considers the case of an elliptic
equilibrium, i.e., the Hamiltonian~\frmref{eq:1}
(see~\dbiref{Whittaker-1916}, \dbiref{Cherry-1924a}, \dbiref{Cherry-1924b} and \dbiref{Birkhoff-1927}).
Extensions to the case of maps have been studied by Poincar\'e, who
introduced the idea of reducing the flow of differential equations to
a map\bibref{Poincare-1890}\bibref{Poincare-1912}.  On the other
hand, Birkhoff exploited the idea of interpolating an area preserving
map of the plane by a Hamiltonian flow\bibref{Birkhoff-1920}.

A wide discussion of the problem of control for Hamiltonian systems in
the light of renormalization theory has been made by
Gallavotti\bibref{Gallavotti-1982}, who calls {\corsivo counterterms}
the control terms.  A constructive approach has been proposed by
Vittot\bibref{Vittot-2004}, by using the Lie series formalism in order
to construct a suitable normal form.  Vittot's ideas have been
exploited in order to find control terms that may reduce chaos (or
increase the size of the stability region around the equilibrium
point) in models of interest in physics, e.g., the dynamics of
magnetized
plasmas\bibref{CirBriChaFloLimVitPet-2004}\bibref{CirChaLimVitPet-2004}\bibref{CirChaLimVitPetFigGhe-2004}\bibref{ChaCirDovLimMacVit-2005}.
Although a full control in the sense initially proposed by Vittot is
unrealistic in a practical application (as we shall discuss below),
introducing some properly chosen terms in the Hamiltonian may
significantly reduce chaos, thus stabilizing the dynamics.  This
aspect, which may present a considerable practical interest, is
discussed in the above quoted papers, with explicit examples.

Concerning maps, the problem has been widely studied since the
eighties of the past century in a series of papers by a group of
authors including Bazzani, Giovannozzi, Servizi, Todesco and
Turchetti.  They investigated the normal form for symplectic maps in
view of application to betatronic motions in
accelerators\bibref{Bazzani-1988}\bibref{Bazzani-1993}\bibref{Bazzani-1994}\bibref{Todesco-1994}.

However, we think that a better insight on the problem may be achieved
by introducing some changes both in the technical tools and in the
theoretical framework.  This is what we are going to discuss in this
paper, making reference in particular to the case of a (symplectic)
map in the neighborhood of an elliptic equilibrium.

The paper is organized as follows.  In the rest of the present section
we include an informal discussion of the technical improvement and of
the theoretical framework.  In section~\secref{sec:formal} we present
the formal algorithm that allows us to calculate the normal form of
the map and of different forms of the control terms.  In
section~\secref{sec:stime} we work out all the quantitative estimates
on the normal form and on the possible control terms.  The numerical
application is presented in section~\secref{sec:numerica}, where we
illustrate an heuristic method to predict the size of the stable
region using the normal form and compare the results with a direct
iteration of the map.

\subsection{sbs:techimpr}{Technical improvements}
The construction of a normal form for symplectic maps is usually
worked out by using interpolation via a canonical system of
differential equations, as proposed by Birkhoff.  The reason is
perhaps that for the Hamiltonian case there are several methods
available, the most effective ones being based on performing canonical
transformation using the method of Lie series, as proposed in the recent paper
of Vittot\bibref{Vittot-2004}.

Let us start by making two remarks.  The first one is that going through a flow of
differential equations in order to represent a map is a lengthy
procedure, that may be desirable to avoid.  The second remark is that
replacing the method of Lie series with that of Lie transforms makes
the calculations (to be done via algebraic manipulation) definitely
more effective and introduces a lot of flexibility.

We adopt here a language that is common in the milieu of Celestial
Mechanics.  Given a vector field $X(x)$ on a $n$-dimensional manifold
with coordinates $x$ the Lie series is introduced as the differential
operator $\exp\bigl(t\lie{X}\bigr)=\sum_{k\ge
0}\frac{t^k\lie{X}^k}{k!}$, where $\lie{X}$ is the Lie derivative
along the flow generated by the vector field $X$ and $t$ is the time
parameter of the flow.  In this form the Lie series has been widely
used by Gr\"obner\bibref{Groebner-1957}\bibref{Groebner-1967}.  The
Lie transform is a generalization of Lie series that has been
introduced in different ways.  We hereby use the following definition.  Let
$X=\{X_j\}_{j\ge 1}$ be a sequence of vector fields, which are
supposed to be of increasing order in some small parameter.  The Lie
transform operator $T_X$ is defined as $T_{X} = \sum_{s\ge 0} E_s$,
where the sequence $E_s$ of linear operators is recursively defined as
$E_0 = \uno$ and $E_s = \sum_{j=1}^{s} \frac{j}{s}\lie{X_j} E_{s-j}$.
By letting the sequence to have only one vector field different from
zero, e.g., $X=\{0,\ldots,0,X_k,0,\ldots\}$ it is easily seen that one
gets $T_X=\exp\bigl(\lie{X_k}\bigr)$.  Thus, a Lie series may be seen
as a particular case of a Lie transform (see sect.~\secref{sec:formal}
below).  Actually the Lie series in the form above, namely with a
homogeneous polynomial vector field, is unable to reproduce any near
the identity transformation for which a composition of Lie series is
necessary.  This problem is overcome with the algorithms introduced by
Hori\bibref{Hori-1966} and Deprit\bibref{Deprit-1969}.  Hori allowed
the vector field $X(x,\epsilon)$ to depend also on a small
perturbation parameter $\epsilon$ and looked for a series expansion of
the flow up to time $\epsilon$, thus introducing the operator
$\exp\bigl(\lie{X(x,\epsilon)}\bigr)$.  A method similar to that of
Hori has been introduced also by Takens\bibref{Takens-1974} and
applied to maps in~\dbiref{GelGel-2009}.  Deprit\bibref{Deprit-1969}
considered a non-autonomous vector field $X(x,t)$, letting the time of
the flow to play the role of a small parameter.  The algorithms so
found are equivalent, although the generating vector field is not the
same in the two cases.  Both of them are able to reproduce any near
the identity transformation. Later, a purely algebraic procedure that
leads to introducing the algorithm for the Lie transform adopted in
the present paper has been proposed in~\dbiref{GioGal-1978} and has
been given a quantitative form in~\dbiref{Giorgilli-1985}.  Actually,
different algorithms representing a Lie transform have been proposed
in the literature: see, e.g.,~\dbiref{HenRoe-1974}.

Our first point is that a holomorphic map in a neighborhood of an
equilibrium may be represented as a composition of a near the identity
transformation with a linear one.  Now, the near the identity
transformation may be represented by a Lie transform, as we have said.
On the other hand, a linear transformation may be seen as the
time-one flow of a linear system of differential equations, which in
turn is represented as a Lie series.  Therefore a map may be
represented as the composition of a Lie transform with a Lie series.
This elementary remark is the key that allows us to get rid of the
interpolation via a flow.  The transformation to normal form is worked
out using a second Lie transform.

A crucial technical tool in this connection is concerned with the
composition of Lie transforms, as opposed to composition of Lie
series.  The point is that the composition of Lie series is 
handled through the well known Baker-Campbell-Hausdorff (BCH)
formula.  I.e., given two linear operators $A$ and $B$ one looks for a
linear operator $C$ such that $e^A\circ e^B= e^C$.  The problem is
that $C$ has a cumbersome form, which makes it difficult to implement
an iterative procedure because it is a hard task to separate terms of 
different orders entering a power expansion in either the
coordinates or a parameter.  In contrast, the composition of two Lie
transforms may be expressed as a third Lie transform in a form that
can be effectively iterated, because it gives an explicit expansion
order by order.  The latter formula is given later in
sect.~\secref{sec:formal}, where also the method of representation of
maps is recalled.  For a short but essentially complete exposition of
the method we refer to~\dbiref{Giorgilli-2013}.

\subsection{sbs:teoria}{Theoretical framework}
Most of the recent works on the theory of control basically attempt to
introduce corrections that reduce the current system to an integrable
one.  As a matter of fact, however, the actual corrections that are introduced
in the examples, having in mind possible physical applications,
consist merely in changing some terms, because a complete correction
may be devised theoretically in many ways, but it can hardly be done
in practice.

Our remark is that integrability is a too strong request: even from
the theoretical viewpoint a better approach can be made based on the
theory of exponential stability developed by Littlewood and Moser, and
fully stated by Nekhoroshev.  Our proposal is to base the search for a
control on stability over long times, i.e., effective stability.

In rough terms, the crucial point of Nekhoroshev's theory is that if
one considers only the actions, $p$ say, of a Hamiltonian system
like~\frmref{eq:2} then one may prove an inequality such as
$\bigl|p(t)-p(0)\bigr|\lt\epsilon^{b}$ for $|t|\lt\exp(1/\epsilon^a)$
with some positive constants $a,b\lt 1$.  The exponential increase of
the time with the inverse of the perturbation is the relevant
phenomenon, because if the perturbation is made small enough then the
estimated time may exceed the lifetime of the physical system that we
are considering.  The proof of the theorem for the general
system~\frmref{eq:2} is not constructive even if one uses the Lie
transform methods, due to the need of a clever geometrical
construction of resonant and non-resonant domains.  However, a simpler
constructive proof can be worked out if one considers the case of an
elliptic equilibrium for a Hamiltonian system, namely a system
like~\frmref{eq:1}.  The same remark applies to the case of symplectic
maps in a neighborhood of an elliptic equilibrium, as we do in the
present paper.

The remarkable advantage of Nekhoroshev's exponential stability is
that it applies to non-integrable systems, and may be very useful if
one is concerned only with long time stability of the motion, while
the actual orbits are not relevant.  Moreover, it applies also if a
small not periodic (smooth) time dependence is added to the
perturbation (see~\dbiref{GioZeh-1992}).  In a practical application this
may be considered to mean that very small perturbations not taken into
account by the model should not be harmful, e.g., slow but very small changes in the
behavior of the magnets of an accelerator.  A
well known problem, that shows up also with Kolmogorov's theory, is
that the analytic estimates on the stability time usually give
unrealistic and even ridiculous results.  However, realistic
estimates may be obtained by complementing the analytic theory with
algebraic manipulation.  This has been verified in the last two
decades through some applications, mainly in the field of Celestial
Mechanics, for both Kolmogorov's
theorem\bibref{Celletti-1988}\bibref{Celletti-1997}\bibref{LocGio-2000}\bibref{LocGio-2007}
and Nekhoroshev's
theory\bibref{GioSko-1997}\bibref{GioLocSan-2009}\bibref{SanLocGio-2010}.

Better estimates are obtained by explicitly calculating the normal
form up to an high order $r$ via algebraic manipulation and by
numerically evaluating the size of the remainder.  Both these
operations are made effective by the use of the Lie transform
algorithm, because the contributions of every order are easily
separated.  The drawback is that also this method produces too
pessimistic results, thus we resort to an heuristic criterion that will be illustrated in section~\secref{sec:numerica}.

\subsection{sbs:add_contro}{Adding a control}
We come now to the formulation of the control theory via normal form.
We refer again to the case of a Hamiltonian system, so that a
comparison with the previous literature is easily done.  However we
stress once more that there is no substantial difference in the case
of maps.

Suppose that we have constructed the normal form 
$$
H^{(r)}(x,y) = H_0 + Z^{(r)}(I_1,\ldots,I_n) + Q^{(r+1)}(x,y)
$$ 
at finite order $r$, as above.  Theoretically, the simplest action is
to choose the control term as $-Q^{(r+1)}$, so that the system turns
out to be integrable.  This is indeed the suggestion made by Vittot,
but his construction is made only for $r=1$.  An attempt to iterate
such a scheme has been made in~\dbiref{VitChaCirLim-2005}. We
emphasize that this may be done in many ways, depending on the choice
of the normalization order $r$, which is the advantage offered by the
availability of a scheme that can be easily iterated.  That is, a
system may be made integrable by adding control terms of order higher
than the first one.  The actual applicability in an experimental
framework remains a major problem, of course, mainly because the
control term provided by the theory should be made using real magnets
(see~\dbiref{BorCarSkoPapVit-2012} for a partial positive answer).

A different but possibly better approach is to introduce some
correction at a given order $s\lt r$ (e.g., the first perturbation
order, as in most of the existing models) and then evaluate the size
of the domain where a long time stability is expected to hold.  Here,
the choice of $s$ is arbitrary, while $r$ should be taken large enough,
compatibly with the computer power.  Successive refinements of the
corrections at different orders are not excluded: this may be easy to
do numerically, but, again, the actual applicability in a physical experiment
may be questionable.

We focus in particular on the use of the normal form in order to
predict the effectiveness of a correction.  The main point is that the
size of the stability region around an equilibrium may be evaluated by
looking at the numerically calculated normal form.  Examples will be
given in sect.~\secref{sec:numerica}.

We add just a remark.  It may seem that implementing a numerical
procedure of calculation of the normal form and repeating it by trying
different corrections should be time consuming, or even impractical.
However, our previous experience shows that this is untrue.  For low
dimensional maps (i.e., in a phase space with a number of degrees of
freedom up to 3) and up to not too high orders the construction of the
normal form with a suitably devised program of algebraic manipulation
it is matter of minutes on a standard personal computer.  The usual method of checking the effect of
corrections by graphical methods (i.e., by just looking at the diagram
of different orbits or by frequency analysis) may require a definitely
longer time, if one wants a complete statistics of orbits.

\subsection{sbs:teorema}{Formal statement}
We give the formal statement of a theorem which represents the basis
of the control method proposed in this paper.  We consider a map of
the form~\frmref{map:1} assuming a non-resonance condition on the
frequency vector $\omega$, namely that $e^{i\langle
k,\omega\rangle \pm i \omega_j}-1\neq0$ for $k \in \interi^n$, $|k|>1$
and for $j=1,\ldots,n$.  Observe that the above are the minimal
assumptions to have a formal statement.
\theorem{thm:teorema}
Let the symplectic map
$$
z' = \Lambdamatw z + F_2(z)\ , \quad z=(x,y)\in\reali^{2n}
$$
with non-resonant frequencies $\omega$ be given, where $F_2(z)$ is a power
series starting with terms of degree at least $2$.  Then for every $r\geq1$ there exists a
formal near the identity symplectic transformation which gives the map the normal form
$$
z' = \Omegamat^{(r)} z + P^{(r+1)}(z)\ ,
\formula{eq:mappanorm}
$$
where
$\Omegamat^{(r)}=\Lambdamatw+\tilde\Omegamat^{(r)}(I_1,\ldots,I_n)$ is
a symplectic matrix with $\tilde\Omegamat^{(r)}(I_1,\ldots,I_n)$ at
least linear in the actions $I_1=\frac{x_1^2+y_1^2}{2}, \ldots,
I_n=\frac{x_n^2+y_n^2}{2}$, and $P^{(r+1)}$ is a power series starting
with terms of degree at least $r+2$.  Finally the truncated
 map $z'=\Omegamat^{(r)} z$ is integrable.
\endclaim

The formal statement of theorem~\thrref{thm:teorema} is sufficient to
have a control on the map.  However, we should stress that our method
is inspired by the long-time stability theory in Nekhoroshev sense.
Therefore we discuss how such an estimate could be worked out by
adding quantitative estimates to the theorem.

A rigorous estimate requires a stronger non-resonance condition that
we take to be the Diophantine one, namely
$$
|e^{i\langle k,\omega\rangle \pm i \omega_j}-1|
 \geq \frac{\gamma}{|k|^\tau}\quad
\hbox{with}\quad\gamma>0\ ,\  \tau>n-1\ .
$$

The suggestion is to transport to the case of maps the known results
concerning Hamiltonian flows.  Precisely, referring to the
Hamiltonian~\frmref{eq:1} constructing the well known Birkhoff normal
form up to a finite order $r>1$ means that the Hamiltonian can be
transformed to
$$
Z^{(r)} = H_0 + Z_1 + \ldots + Z_r + Q^{(r+1)}\ ,
$$
where $Z_s(I_1,\ldots,I_n)$ depends only on the actions
$I_1=\frac{x_1^2+y_1^2}{2}, \ldots, I_n=\frac{x_n^2+y_n^2}{2}$.
Formally this corresponds to theorem~\thrref{thm:teorema} above.
Assume that the frequencies are Diophantine, i.e., $\bigl|\langle
k,\omega\rangle\bigr|\geq\frac{\gamma}{|k|^\tau}$, with $\gamma$ and
$\tau$ as above.  Then a quantitative version states that in a
neighborhood of the origin of radius $\rho$ the estimate
$$
\bigl|Q^{(r+1)}\bigr|_\rho < C^r r!^{\tau+1} \rho^{r+3}
$$
holds true for the supremum norm of the remainder $Q^{(r+1)}$.

The corresponding result for a map may be shortly stated as
\proposition{pro:quantitativo}
There exists $\rho^*(r)$ such that the transformation to normal form
is holomorphic in a neighborhood of radius $\rho<\rho^*(r)$ around
the origin and the norm of the remainder $P^{(r+1)}(z)$
in~\frmref{eq:mappanorm} is estimated as
$$
\bigl|P^{(r+1)}(z)\bigr| < C^r r!^{\tau+1} \rho^{r+2}
\formula{frm:stimaresto}
$$
with some constant $C>0$.
\endclaim

\noindent
The relevance of the proposition is that the
estimate~\frmref{frm:stimaresto} of the not normalized remainder has
precisely the form that allows us to get exponential stability
estimates, as sketched in sect.~\sbsref{sbs:teoria}.  Indeed, for
every radius $\rho$, one can select an optimal value $r_{\rm opt}$ of
the normalization order so as to minimize the remainder in
~\frmref{frm:stimaresto}.  Thus one gets an exponential estimate.  We
omit a formal statement of the latter part, since we are mainly
interested in producing a computer assisted method.

\section{sec:formal}{Formal setting}
This section is devoted to a discussion of some formal methods of
perturbation theory that will be used in the rest of the paper.

\subsection{sbs:lie}{Short reminder about Lie series and Lie transform methods}
We briefly recall the definitions of Lie series and Lie transform,
restricting our attention to the case of polynomial vector fields.
Here we include only the essential information for this paper,
referring for a more complete exposition to~\dbiref{Groebner-1967}
and~\dbiref{Giorgilli-2013}.  We consider in particular the case of
vector fields and the application to maps.

Let $X_s(z)$ be a vector field on $\complessi^n$ whose components are
homogeneous polynomials of degree $s+1$.  We shall say that $X_s(z)$
is of {\corsivo order} $s$, as indicated by the label.  Moreover, in
the following we shall denote by $X_{s,j}$ the $j$-th component of the
vector field $X_s\,$.  The {\corsivo Lie series} operator is defined
as
$$
\exp(\lie{X_s}) = \sum_{j\ge 0} \frac{1}{j!} \lie{X_s}^j
\formula{sielie.100}
$$
where $\lie{X_s}$ is the Lie derivative with respect to the vector
field $X_s$.

Let now $X=\{X_j\}_{j\ge 1}$ be a sequence of polynomial vector fields
of degree $j+1$.  The {\corsivo Lie transform} operator $T_X$ is
defined as
$$
T_{X} = \sum_{s\ge 0} E^{X}_s\ ,
\formula{trslie.1}
$$
where the sequence $E^{X}_s$ of linear operators is recursively
defined as
$$
E^{X}_0 = \uno\ ,\quad 
E^{X}_s = \sum_{j=1}^{s} \frac{j}{s}\lie{X_j} E^{X}_{s-j}\ .
\formula{trslie.2}
$$
The superscript in $E^{X}$ is introduced in order to specify which
sequence of vector fields is intended.  By letting the sequence to
have only one vector field different from zero, e.g.,
$X=\{0,\ldots,0,X_k,0,\ldots\}$ it is easily seen that one gets
$T_X=\exp\bigl(\lie{X_k}\bigr)$.

The definitions above are transported in a straightforward manner to
the symplectic (or Hamiltonian) framework.  In this case we use the
notation $z=(x,y)\in\reali^{2n}$ and use the vector field
$X=\bigl(\parder{\chi}{y},-\parder{\chi}{x}\bigr)$ where $\chi(x,y)$
is a homogeneous polynomial defined on the phase space $\reali^{2n}$
(or, possibly, $\complessi^{2n}$) with canonical coordinates $(x,y)$.

The Lie series and Lie transform are linear operators acting on the
space of holomorphic functions and of holomorphic vector fields.  They
preserve products between functions and commutators between vector
fields, i.e., if $f,\,g$ are functions and $v,\,w$ are vector fields
then one has
$$
T_X (fg)
=T_X f\cdot T_X g\ ,\quad
T_X \{v,w\} = \bigl\{T_X v, T_X w\bigr\}\ ,
\formula{lser.6}
$$
where $\{\cdot,\cdot\}$ is the commutator between vector fields.  If
the vector fields $v$, $w$ are generated by the Hamiltonian functions
$V$, $W$, then the commutator $\{v,w\}$ is the vector field generated
by the Poisson bracket $\{V,W\}$.  Here, replacing $T_X$ with
$\exp\bigl(\lie{X}\bigr)$ gives the corresponding properties for Lie
series.  Moreover both operators are invertible.

We recall a remarkable property which justifies the usefulness of
Lie methods in perturbation theory.  We adopt the name {\corsivo
exchange theorem} introduced by Gr\"obner.  Let $f$ be a function and
$v$ be a vector field. {\corsivo Consider the 
transformation $w = T_X z$,
i.e., in coordinates,
$$
w_j = z_j + X_{1,j}(z) 
 +  \left[\frac{1}{2} \lie{X_1} X_{1,j}(z) + X_{2,j}(z)\right] 
  +\ldots
\ ,\quad j=1,\ldots,n\ .
$$
and denote by $\Jmat$ its differential, namely, in
coordinates, the Jacobian matrix with elements $J_{j,k}
=\derpar{w_{j}}{z_k}$.  Then one has
$$
f(w)\big|_{w=T_X z} = 
 \bigl(T_X f\bigr) (z)\ ,\quad
\Jmat^{-1} v(w)\Big|_{w= T_X z} = 
 \bigl(T_X v\bigr)(z)\ .
\formula{trslie.4}
$$}

\noindent
That is: there is no need to perform a substitution of variables
as required by the left member; just transform the function of vector
field, as in the right member.

The same statement holds true for a transformation generated by a Lie
series.

\subsection{sbs:rappresentazione}{Representation of maps}
We consider a nonlinear analytic map in the neighborhood of a fixed
point, that we may assume to be the origin.  More precisely we write
the map in general form as
$$
z' = \Lambdamatw z + \phi_1(z) + \phi_2(z) + \ldots\ ,
$$
where $z\in\complessi^n$ and $\phi_s(z)$ is a homogeneous polynomial
of degree $s+1$.  For problems related to control, it is more natural
to consider a symplectic map in canonical variables
$z=(x,y)\in\reali^{2n}$ (or possibly $\complessi^{2n}$) where
$\Lambdamatw$ is a rotation matrix.  A typical form of such a map is
$$
\pmatrix{x' \cr y'} =
\Lambdamatw
\pmatrix{x + f_1(x,y) + f_2(x,y) + \ldots \cr y + g_1(x,y) + g_2(x,y) + \ldots}
\formula{map.1}
$$
with the block matrix
$$
\Lambdamatw = 
\pmatrix{
C_\omega & -S_\omega\cr
S_\omega & C_\omega\cr
}\ .
\formula{eq:lambda}
$$
where $C_\omega=\diag(\cos\omega_1,\ldots,\cos\omega_n)$ and
$S_\omega=\diag(\sin\omega_1,\ldots,\sin\omega_n)$. Usually the
functions $f_j$ and $g_k$ must satisfy the known symplecticity
condition, although this is not strictly necessary for most of our
discussion.  As a typical example one may consider the well known
quadratic map of H\'enon
$$
\pmatrix{x' \cr y'} =
\pmatrix{
\cos\omega & -\sin\omega \cr
\sin\omega &  \phantom{-}\cos\omega\cr}
\pmatrix{x \cr y -x^2}\ .
$$

We recall now the representation of maps introduced
in~\dbiref{Giorgilli-2013} together with some formal results that we
are going to use here.  A trivial but useful remark is the following.
Let $\Lambdamatw=e^{\Amat_\omega}$.  Then we may express the
linear part of the map as a Lie series by introducing the exponential
operator $\Rmatw=\exp\bigl(\lie{\Amat_\omega z}\bigr)$.  The action of the
operator $\Rmatw$ on a function $f$ or on a vector field $V$ is easily
calculated as
$$
\bigl(\Rmatw f\bigr)(z) = f(\Lambdamatw z)\ ,\quad
 \bigl(\Rmatw V\bigr)(z) = \Lambdamatw^{-1} V(\Lambdamatw z)\ ,
\formula{sielie.112}
$$
namely by direct substitution.

The first result is concerned with the representation of the
map~\frmref{map.1} using a Lie transform.

\lemma{map.13}{There exist generating sequences of
vector fields $V= \bigl\{V_s(z)\bigr\}_{s\geq1}$ and $W
=\bigl\{W_s(z)\bigr\}_{s\geq1}$ with $W_s =\Rmatw V_s$ such that
the map~\frmref{map.1} is represented in either form
$$
z' = \Rmatw\circ T_{V} z\quad \hbox{or}\quad
z' = T_{W}\circ\Rmatw\, z
\formula{map.5}
$$
}\endclaim

\noindent
A similar representation, namely the composition of a rotation with a
flow expressed as Lie series, has been introduced by
Gelfreich~et~al. in~\dbiref{GelGel-2009} in order to represent an
area-preserving map of the plane around an elliptic point. From the
theoretical point of view our method is equivalent to that presented
in~\dbiref{GelGel-2009}.  However since they use the Lie series with a
vector field expanded itself in a series, it seems to be less adapted
to an explicit calculation, since a natural identification of
homogeneous terms is not straightforward.

The second result is concerned with the composition of Lie
transforms. 

\lemma{trslie.11}{Let $X,\,Y$ be generating sequences.  Then one
has $T_X\circ T_Y = T_Z$ where $Z$ is the generating sequence
recursively defined as
$$
Z_1 = X_1 + Y_1\ ,\quad
Z_s = X_s + Y_s + \sum_{j=1}^{s-1} \frac{j}{s} E^{X}_{s-j} Y_j\ .
\formula{trslie.12}
$$
}\endclaim

\noindent
The latter formula reminds the well known Baker-Campbell-Hausdorff
composition of exponentials.  The difference is that the result is
expressed as a Lie transform instead of an exponential, which makes
the formula more effective for our purposes, as remarked in the
introduction.

The construction of the normal form is based on conjugation of maps.
The question may be stated as follows.  Let two maps
$$
w'=T_{W}\circ\Rmatw\, w\ ,\quad z'=T_{Z}\circ \Rmatw\, z
\formula{map.6}
$$ 
be given, where $\Rmatw$ is a Lie series operator and
$W=\{W_1,W_2,\ldots\}$, $Z=\{Z_1,Z_2,\ldots\}$  are generating
sequences.  The problem is to find
whether the maps are conjugated by  a near the identity transformation
$$
w = z + \phi_1(z) + \phi_2(z) +\ldots\ .
\formula{map.7}
$$
Using the same operator $\Rmatw$ in both maps means only that the
unperturbed maps are trivially conjugated.  The question may be
answered either in formal sense (roughly: disregarding the convergence
of the series) or asking the stronger property that the transformation
is holomorphic.

\lemma{map.10}{Let  $X=\{X_1,X_2,\ldots\}$ be a generating
sequence of the near the identity transformation $w=T_X z$.  Then the
maps~\frmref{map.6} are formally conjugated if 
$$
T_{Z}\circ T_{\Rmatw X} = T_{X} \circ T_{W}\ .
\formula{map.9}
$$
More explicitly, the following relations must be satisfied:
\formdef{map.20c}
\formdef{map.20b}
$$
\displaylines{
\frmref{map.20c}\hfill
\Dmatw X_1 + Z_1 = W_1\ ,\quad \Dmatw = \Rmatw-\uno
\hfill\cr
\frmref{map.20b}\hfill
\Dmatw X_s + Z_s=
 W_s + \sum_{j=1}^{s-1}\frac{j}{s} 
  \bigl(E^{X}_{s-j} W_j - E^{Z}_{s-j} \Rmatw X_j\bigr)  
\ ,\quad s\gt 1\ .
\hfill\cr
}
$$}\endclaim

\subsection{sbs:nfmappe}{Normal form and control for maps}
A normal form for the map may be constructed as follows.
We use the conjugation formula of lemma~\lemref{map.10} by asking the
generating sequence $Z$ to possess some nice property that we shall
specify later.  In general, this process leads to divergent series,
but we can stop the construction of the normal form by picking $r\geq1$
and looking for a finite generating sequence
$X^{(r)}=\{X_1,\ldots,X_r,0,\ldots\}$ such that the transformed vector
field $Z$ is in normal form up to order $r$.  That is we want
$Z=\{Z_1,\ldots,Z_r,Q_{r+1},\ldots\}$ where
$\Qscr=\{0,\ldots,0,Q_{r+1},Q_{r+2},\ldots\}$ is a not normalized
remainder.

The common way to introduce a control is to look for a map
generated by a vector field $V=W+\Fscr^{(r+1)}$ where
$\Fscr^{(r+1)}=\{0,\ldots,0,F_{r+1},F_{r+2},\ldots\}$, such that $V$ and
$W$ coincide up to order $r$.  The control $\Fscr^{(r+1)}$ is
determined so that the normal form for the controlled vector field $V$
is a finite sequence $Z^{(r)}=\{Z_1,Z_2,\ldots,Z_r,0,\ldots\}$, i.e., $Q_{r+1} =
Q_{r+2} = \ldots = 0$.

Assuming that $W$ is known, the normal form $Z$ is determined by
solving for $Z_1,\ldots,Z_r$ and $X_1,\ldots,X_r$ the
equations~\frmref{map.20c} and~\frmref{map.20b} of
lemma~\lemref{map.10}, for $s=1,\ldots,r\,$.

Then the control term
$\Fscr^{(r+1)}$ is determined as
$$
F^{(r+1)}_s = -W_s - \sum_{j=1}^{s-1} \frac{j}{s}
\bigl(E^{X^{(r)}}_{s-j} \tilde{W}_j - E^{Z^{(r)}}_{s-j} \Rmatw X^{(r)}_j\bigr)  
\ ,\quad s\gt r\ ,
\formula{eq:controllo}
$$
where $\tilde{W}_j = {W}_j$ for $j\leq r$ and $\tilde{W}_j =
{W}_j + F^{(r+1)}_j$ for $j> r\,$.

We remark that this is indeed an extension of the method used by
Vittot\bibref{Vittot-2004} which, however, is based on Lie series.
Actually, the method of Vittot is equivalent to just introducing a
control $\Fscr^{(2)}=\{0,F_{2},0,\ldots\}$. This was given, e.g.,
in~\dbiref{BorCarSkoPapVit-2012}.  Let us observe that the use of Lie
transforms is more suitable than the Lie series, to get high order
control terms.  The correspondence is not straightforward, as it
happens also for the correspondence between the Lie series and Lie
transform (see, e.g., \dbiref{Fasso-1989}).

It remains to discuss how to solve a homological equation of the form
$$
\Dmatw X + Z = \Psi\ ,\quad \Dmatw = \Rmatw-\uno
\formula{eq:homeq-map}
$$
where $\Psi$ is a known homogeneous polynomial of degree $s\geq2$ and
$X$ and $Z$ are the unknowns.

Referring to a system of the form~\frmref{map.1}, the main remark is
that the matrix $\Lambdamatw$ takes a diagonal form via the
transformation
$$
x_l =\frac{1}{\sqrt{2}}(\xi_l +i\eta_l) \ , \quad y_l
=\frac{i}{\sqrt{2}}(\xi_l- i\eta_l) \ , \quad 1\le l\le n \ .
\formula{eq:complex}
$$
We remark that the transformation is symplectic.  Then
$$
\Lambdamatw=
\diag(e^{-i\omega_1},\ldots,e^{-i\omega_n}, e^{i\omega_1},\ldots,e^{i\omega_n})
$$
turns out to be diagonal.  Thus, the
operator $\Dmatw$ in~\frmref{eq:homeq-map} turns out to be diagonal,
too.  This is seen as follows.  Denote by
$(\evet_1,\ldots,\evet_{2n})$ the canonical basis of
$\complessi^{2n}$, and expand the vector field $\Psi$ over the basis
of monomials $\xi^j \eta^k \evet_l$.  Then a straightforward
calculation gives
$$
\Dmatw \xi^j \eta^k \evet_l =
(e^{i\langle(k-j),\omega\rangle + i \mu_l}-1)\, \xi^j \eta^k \evet_l\ ,
\formula{eq:mu}
$$
where $\mu = (-\omega_1,\ldots,-\omega_n,\omega_1,\ldots,\omega_n)$.
Here we used the notation $\langle k,\omega\rangle=\sum_{l=1}^{n}
k_l \omega_l\,$. Thus, dealing with the homological equation is a
trivial matter.

\subsection{sbs:provateoremaformale}{Proof of theorem~\thrref{thm:teorema}}
Using lemma~\lemref{map.13}, we can represent the map as $z' =
T_{W}\circ\Rmatw\, z$.  Using lemma~\lemref{map.10}, we look for a
truncated generating sequence $X^{(r)}=\{X_1,\ldots,X_r,0,\ldots\}$
that conjugates the map to $z' = T_{Z} z$ with
$Z=\{Z_1,\ldots,Z_r,Q_{r+1},\ldots\}$, with the further request that
$Z_1,\ldots,Z_r$ are in Birkhoff normal form.  That is they commute
with $\Rmatw$, or equivalently, they depend only on the actions
$I_1,\ldots,I_n$ in view of the non-resonance condition on $\omega$.
We truncate the sequence $Z$ by defining $\Gamma^{(r)}
= \{Z_1,\ldots,Z_r,0,\ldots\}$ and consider the map
$$
\tilde z = T_{\Gamma^{(r)}}\circ\Rmatw\, z\ ,
$$
which is clearly symplectic.  On the other hand, it is immediately
checked that the latter map has the form
$$
\tilde z = \Omegamat^{(r)} z\ ,
$$
as claimed, so that it is obviously integrable.  Thus we recover the
transformed map
$$
z' = \tilde z + P^{(r+1)}(z)\ ,
$$
with $P^{(r+1)}(z)=T_{Z}\circ\Rmatw\, z -
T_{\Gamma^{(r)}}\circ\Rmatw\, z$ which clearly is a series starting
with terms of degree at least $r+2$.  This concludes the proof.

\subsection{sbs:hom}{One possible solution of the homological equation}
We come now to discuss how to solve the homological
equation~\frmref{eq:homeq-map}.  We recall that $\Psi$ is a
homogeneous polynomial of some degree $r\geq2$ and $X$, $Z$ are the
unknown vector fields that we want to determine as homogeneous
polynomial of the same degree.

The traditional approach is the following. Denote by $\Pi_r$ the space
of homogeneous polynomial vector fields of fixed degree $r$.  The
kernel $\Nscr$ and the range $\Rscr$ of the linear operator $\Dmatw$
are defined as usual, namely
$$
\Nscr = \Dmatw^{-1} 0\ ,\quad
\Rscr = \Dmatw \Pi_r\ .
$$
Since $\Dmatw$ maps $\Pi_r$ into itself and can be diagonalized, the
properties $\Nscr\cap\Rscr=\{0\}$ and $\Nscr\oplus\Rscr=\Pi_r$ hold
true.  Thus the inverse $\Dmatw^{-1}$ is well defined if we consider
$\Dmatw$ restricted to the range $\Rscr$.  The projectors
$\projn$ and $\projr$ are naturally introduced as
$$
\projr = \Dmatw^{-1}\circ\Dmatw \ ,\quad
\projn = \uno-\projr\ .
$$

The range and the kernel are explicitly characterized as follows.
Define the resonance module
$$
\Mscr_\omega = \{ k\in\interi^{n}\>:\> e^{i \langle k,\omega \rangle}-1 =
0  \}\ .
$$
Recall that a symplectic vector field is generated as
$X_{H}=\Jmat \nabla H$ where $H$ is a Hamiltonian function.  Consider
now a monomial $H=\xi^j \eta^k$.  Then $X_{H}=(U, V)^{T}$ with
$$
U_l = -k_l \xi^j\eta^{k-d_l}
\quad\hbox{and}\quad
V_l = j_l \xi^{j-d_l}\eta^{k}\ ,\quad l=1,\ldots,n\ ,
$$
where $d_l$ is a vector of the canonical basis for $n$-dimensional integers.  Therefore
$$
\Dmatw X_H = 
(e^{\langle k-j,\omega\rangle}-1) X_H\ .
$$
The resonant terms in $X_H$ correspond to the monomials
$\xi^j\eta^{k-e_l}$ and $\xi^{j-e_l}\eta^{k}$ with
$k-j\in \Mscr_\omega$.

A straightforward solution of equation~\frmref{eq:homeq-map} is found
by setting $Z=\projn \Psi$ and $X = \Dmatw^{-1}\projr \Psi$.  In other
words, we put into $Z$ all the resonant monomials found in $\Psi$, and
solve the homological equation only for the non-resonant monomials.
This ensures that $\Dmatw Z=0$, which can actually be used as the
condition for $Z$ to be in normal form.  The procedure outlined here
is the natural adaptation of the well known Birkhoff normal form for
Hamiltonian systems to the case of maps.  It is an easy matter to prove
that if the original map is symplectic, so is the normal form.

As a matter of fact the solution above is very restrictive.  Other
characterizations may be devised, depending precisely on how
equation~\frmref{eq:homeq-map} is solved.

\section{sec:stime}{Quantitative estimates}
In this section we produce quantitative estimates for the generating
sequence of the normal form.  The aim is to provide to the reader
enough information in order to work out a proof of
proposition~\proref{pro:quantitativo}.  Actually the more difficult
part is to give estimates for the generating sequence.  The final part
of the proof is just an adaptation of the proofs already published in
previous papers, see,
e.g.,~\dbiref{Giorgilli-1985}, \dbiref{GioDelFinGalSim-1989}
or \dbiref{Giorgilli-2003}.  Hence for this part we give only a short
sketch.

\subsection{sbs:divisori}{Small divisors}
Controlling the divisors in the solution of the homological equation
is usually a major task that can be hardly worked out in a general
way.  For instance, which divisors do occur and how they do accumulate depends
on the particular normal form that one is looking for.
The estimates below are particularly adapted to the case of Birkhoff
normal form, both in the resonant and in the non-resonant case.

Let us refer again to the map~\frmref{map.1}.  For a given vector
$\omega\in\reali^n$ consider the set of non-negative integers vectors
and labels which are non-resonant with $\omega$.  Formally:
$$
\Kscr_{\omega} = \bigl\{(k,j)\in\interi_+^{n}\times \{1,\ldots,n\}\>:\> 
 e^{i\langle k,\omega\rangle \pm i \omega_j}-1 \ne 0  \bigr\}
\formula{map.20}
$$
Let the positive sequences $\{\beta_s\}_{s\gt 0}$ and $\{\alpha_s\}_{s\ge 0}$ be
defined as
$$
\beta_s = \min_{{|k|=s+1}\atop{(k,j)\in\Kscr_{\omega}}} 
 \bigl|e^{i\langle k,\omega\rangle \pm i \omega_j}-1\bigr|
\formula{map.21}
$$
and
$$
\alpha_0 = 1\ ,\quad
 \alpha_s = \min (\beta_s,\alpha_{s-1})\ . 
\formula{map.22}
$$
Thus, $\alpha_s$ is the smallest divisor that may appear in the
solution of the homological equation up to order $s$.  If the strong
non-resonance is assumed that the frequencies are Diophantine, then we
have
$$
\alpha_s \geq \frac{\gamma}{s^\tau}\quad\hbox{with}\quad\gamma>0\ ,\  \tau>n-1\ .
$$

We emphasize that the definitions of the set $\Kscr_{\omega}$ and of
the sequences above strongly depends on the choice made for
characterizing the normal form.

Finally let us define the useful sequence $\{T_{r}\}_{r\ge 0}$ as
$$
T_{0} = 1\ , \quad T_r = \frac{1}{\alpha_r}T_{r-1}
$$
that will be used to control the accumulation of the divisors (see
lemma~\lemref{sielie.68} below).

\subsection{sbs:norme}{Norms on vector fields and generalized Cauchy estimates}
For a homogeneous polynomial $f(z)=\sum_{|k|=s}f_kz^k$ (using
multiindex notation and with $|k|=|k_1|+\ldots+|k_n|$) with complex
coefficients $f_k$ and for a homogeneous polynomial vector field
$X_s=(X_{s,1},\ldots,X_{s,n})$ we use the {\it polynomial norm}
$$
\bignorma{f} = \sum_{k} |f_k|\ ,\quad
 \bignorma{X_s} = \sum_{j=1}^{n}\, \bignorma{X_{s,j}}\ .
\formula{sielie.104}
$$

The following lemma allows us to control the norms of Lie derivatives
of functions and vector fields.
 
\lemma{sielie.105}{Let $X_r$ be a homogeneous polynomial vector field
of degree $r+1$.  Let $f_s$ and $v_s$ be a homogeneous polynomial and
vector field, respectively, of degree $s+1$.  Then we have
$$
\bigl\|\lie{X_r} f_s\bigr\| \le 
 (s+1)\, \|X_r\|\, \|f_s\|\quad {\rm and}\quad
  \bignorma{\liec{X_r} v_s}
 \le (r+s+2)\,
  \|X_r\|\, \|v_s\|\ .
\formula{sielie.106}
$$}\endclaim

\proof
Write $f_s=\sum_{|k|=s+1} b_k z^k$ with complex coefficients
$b_k$. Similarly, write the $j$-th component of the vector field
$X_r$ as $X_{r,j}=\sum_{|k'|=r+1} c_{j,k'} z^{k'}$.  Recalling that
$\lie{X_r}f_s=\sum_{j=1}^{n}X_{r,j}\derpar{f_s}{z_j}$  we have
$$
\lie{X_r}f_s = \sum_{j=1}^{n} \sum_{k,k'}
 \frac{c_{j,k'}k_jb_k}{z_j} z^{k+k'}\ .
$$
Thus in view of $|k_j|\leq s+1$ we have
$$
\bignorma{\lie{X_r}f_s} 
 \le 
   (s+1)\sum_{j=1}^{n}\sum_{k'} |c_{j,k'}| \sum_{k} |b_k| 
   = (s+1)\ordnorma{X_r}\, \ordnorma{f_s}\ ,
$$
namely the first of~\frmref{sielie.106}.  In order to prove the second
inequality recall that the $j$-th component of the Lie derivative of
the vector field $v_s$ is
$$
\bigl(\lie{X_r} v_s\bigr)_j =
 \sum_{l=1}^{n}\left(X_{r,l}\derpar{v_{s,j}}{z_l} 
  - v_{s,l}\derpar{X_{r,j}}{z_l}\right)\ .
$$
Then using the first of~\frmref{sielie.106} we have
$$
\eqalign{
\biggnorma{\sum_{l=1}^{n}\left(X_{r,l}\derpar{v_{s,j}}{z_l} 
  - v_l\derpar{X_{r,j}}{z_l}\right)}
\le  (s+1) \ordnorma{X_r}\, \ordnorma{v_{s,j}}
    +(r+1) \ordnorma{v_s}\, \ordnorma{X_{r,j}}\ , 
\cr
}
$$
which readily gives the wanted inequality in view of the
definition~\frmref{sielie.104} of the polynomial norm of a vector
field.\endproof

\lemma{sielie.109}{Let $V_r$ be a homogeneous polynomial vector field
of degree $r+1$.  Then the equation $\Dmatw X_r
= \projr V_r$ possesses a unique solution $X_r\in{\cal R}$ satisfying
$$
\bignorma{X_r} \le \frac{1}{\alpha_r}\bignorma{V_r}\ ,\quad
 \bignorma{\Rmat X_r} \le \frac{1+\alpha_r}{\alpha_r}\bignorma{V_r}\ ,
\formula{sielie.110}
$$
with the sequence $\alpha_r$ defined by~\frmref{map.22}
}\endclaim

\proof
The first inequality is a straightforward consequence of the
definition~\frmref{sielie.104} of the norm and of the sequences
$\beta_r$ and $\alpha_r$ defined by~\frmref{map.21}
and~\frmref{map.22}.  For if $v_{j,k}$ are the coefficients of $V_r$
then the coefficients of $X_r$ are bounded by $|v_{j,k}|/\beta_r\le
|v_{j,k}|/\alpha_r$.  The second inequality follows from $\Rmat X_r =
X_r+V_r\,$.\endproof

\subsection{sbs:iteration}{Estimate of the generating sequence}
This is the main lemma which allows us to control the norms of the
sequence of vector fields $\{X_r\}_{r\gt 1}\,$.

\lemma{sielie.68}{Assume that the sequence $W$ of vector fields
satisfies $\ordnorma{W_s}\le C^{s-1}A$ with some
constants $A\gt 0$ and $C\ge 0$.  Then the truncated sequence of vector fields
$X^{(r)}$ that gives the normal form $Z^{(r)}$
satisfies the following estimates:
$$
\bignorma{X_s} \le T_{s} \frac{B_r^{s-1}A}{s}
 \ ,\quad
  \bignorma{Z_s} \le T_{s-1} \frac{B_r^{s-1} A}{s}\ .
\formula{sielie.71}
$$
where
$$
B_r = 4C +8(r+2)A \ .
\formula{sielie.72}
$$
}\endclaim

\proof
Denoting by $\Psi_s$ the r.h.s.\ of~\frmref{map.20c}
and~\frmref{map.20b}, we look for three sequences $\{\eta_s\}_{1\leq
s \leq r}$, $\{\tilde\theta_{s-j,j}\}_{1\leq j\leq s \leq r}$
and $\{\tilde\xi_{s-j,j}\}_{0\leq j \leq s \leq r}$ such that
$$
\eqalign{
&\|\Psi_s\| \leq \eta_s A\, T_{s-1}\ ,\cr
&\|E_{s-j}^{(X)}W_j\| \leq \tilde\theta_{s-j,j} A\, T_{s-j}\ ,\cr
&\|E_{s-j}^{(Z)}\Rmat X_j\| \leq \tilde\xi_{s-j,j} A\, T_{s-1}\ .
}
\formula{eq:stimelemma}
$$
For $s=1$, $\|\Psi_1\|=\|W_1\|\leq A$ and we just set $\eta_1=1$.
Moreover $\|Z_1\|\leq\|\Psi_1\|\leq A$ and, in view of
lemma~\lemref{sielie.109}, $\|X_1\|\leq \|\Psi_1\|/\alpha_1\leq A
T_1$.  Furthermore $\|E_0^{(X)}W_j\| = \|W_j\| \leq \|\Psi_j\|$ and
$\|E_0^{(Z)}\Rmat X_j\| = \|\Rmat
X_j\|\leq \frac{1+\alpha_j}{\alpha_j} \|\Psi_j\|$, thus we set
$\tilde\theta_{0,j}=\eta_j$ and $\tilde\xi_{0,j}=2\eta_j$.

For $s>1$ we have
$$
\eqalign{
\|\Psi_s\| &\leq \|W_s\| + \sum_{j=1}^{s-1} \frac{j}{s}
\bigl(\|E_{s-j}^{(X)}W_j\| + \|E_{s-j}^{(Z)}\Rmat X_j\|\bigr)
\cr
&\leq C^{s-1}A + \sum_{j=1}^{s-1} \frac{j}{s}
(\tilde\theta_{s-j,j}A T_{s-j} + \tilde\xi_{s-j,j} A T_{s-1} )\cr
&\leq \biggl(C^{s-1} + \sum_{j=1}^{s-1} \frac{j}{s}
(\tilde\theta_{s-j,j} + \tilde\xi_{s-j,j})\biggr) A T_{s-1} \ ,
}
$$
where in the last inequality we use the non-decreasing character of
the sequence $T_s$.  Then we recursively define $\eta_s$ as
$$
\eta_s = C^{s-1} + \sum_{j=1}^{s-1} \frac{j}{s}\, (\tilde\theta_{s-j,j} + \tilde\xi_{s-j,j})\ .
\formula{eq:ricorrenti}
$$
Proceeding in a similar way, we write
$$
\eqalign{
\|E_{s-j}^{X} W_j\| &\leq \sum_{l=1}^{s-j} \frac{l}{s-j}
\|\lie{X_l} E^{X}_{s-j-l} W_j\|\cr
&\leq \sum_{l=1}^{s-j} \frac{l}{s-j} (s+2) \|X_l\| \|E^{X}_{s-j-l} W_j\|\cr
&\leq \sum_{l=1}^{s-j} \frac{l}{s-j} (s+2) \frac{\eta_l A T_{l-1}}{\alpha_l} \tilde\theta_{s-j-l,j} A T_{s-j-l}\cr
&\leq \biggl(\frac{s+2}{s-j} \sum_{l=1}^{s-j}  l \eta_l \tilde\theta_{s-j-l,j} A \biggr) A T_{s-j}\ ,
}
$$
and we define
$$
\tilde\theta_{s-j,j} = \frac{s+2}{s-j} \sum_{l=1}^{s-j}  l \eta_l \tilde\theta_{s-j-l,j} A\ .
$$
We also define
$$
\tilde\xi_{s-j,j} = \frac{s+2}{s-j} \sum_{l=1}^{s-j}  l \eta_l \tilde\xi_{s-j-l,j} A\ .
$$
Furthermore, one easily verifies by induction that
$\tilde\theta_{s-j,j} = \eta_j \tilde\theta_{s-j,1}$ and
$\tilde\xi_{s-j,j} = \eta_j \tilde\xi_{s-j,1}$; thus we introduce new
sequences $\theta_{s} = \tilde\theta_{s,1}$ and $\xi_{s}
= \tilde\xi_{s,1}$.  Recalling that
$$
\eta_{s-j} = C^{s-j-1} + \frac{1}{s-j} \sum_{l=1}^{s-j-1} l\eta_l (\theta_{s-j-l} + \xi_{s-j-l})\ ,
$$
one gets the equality
$$
\theta_{s-j} + \xi_{s-j} - (s+2) A \eta_{s-j} =
(s+2)A \eta_{s-j} (\theta_0 + \xi_0) - (s+2)A C^{s-j-1}\ ,
$$
and so the bound
$$
\theta_{s-j} + \xi_{s-j} < 4(s+2)A \eta_{s-j}\ .
$$
Plugging this estimate in~\frmref{eq:ricorrenti}
one  gets\footnote{\dag}{Use
$\sum_{j=1}^{s-1}j \eta_j \eta_{s-j}=\sum_{j=1}^{s-1}(s-j)\eta_j \eta_{s-j}=(s/2)\sum_{j=1}^{s-1}\eta_j \eta_{s-j}\,$.}
$$
\eqalign{
\eta_s &\leq C^{s-1} + 4(s+2)A \frac{1}{s}\sum_{j=1}^{s-1}j \eta_j \eta_{s-j}\cr
&\leq C^{s-1} + 2(s+2)A \sum_{j=1}^{s-1} \eta_j \eta_{s-j}\cr
&\leq (C+2(s+2)A)^{s-1}\mu_s\ ,
}
$$
where
$$
\mu_1=1\ ,\quad
\mu_s = \sum_{j=1}^{s-1}\mu_j\mu_{s-j}\leq \frac{4^{s-1}}{s}\ .
$$
The claim follows from this inequality, remarking that
for $s\leq r$ we can bound the last sum in a geometric way introducing
the constant $B_r$ as in~\frmref{sielie.72}.  It remains  to check
that the sequence $T_s$ controls the accumulation of the small
divisors.  This is a trivial matter that is left to the reader.\endproof

\subsection{sbs:provateorema}{Sketch of the proof of proposition~\proref{pro:quantitativo}}
The proof of the estimate~\frmref{frm:stimaresto} is based on the
convergence of the near the identity transformation $w=T_Y z$, under
the hypothesis that the generating sequence $Y$ satisfies a suitable
convergence condition.  We proceed as follows.  We consider a poly-disk
$\Dscr_{\rho}$ with radius $\rho$ around the origin.  Then the
generating sequence $Y=\{Y_s\}_{s\geq1}$, with $Y_s$ a homogeneous
polynomial vector field of degree $s+1$, is bounded by the norm
$$
|Y_s|_\rho
 = \sup_{z\in\Dscr_{\rho}} \|Y_s(z)\| \leq \|Y_s\| \rho^{s+1}\ .
$$
Then we use the following
\lemma{lem:convergenza}
Consider the generating sequence $Y=\{Y_s\}_{s\geq1}$, with $Y_s$ a
homogeneous polynomial vector field of degree $s+1$.  In a poly-disk
$\Dscr_{\rho}\,$ with radius $\rho$ around the origin assume
$$
|Y_s|_\rho \leq \frac{b^{s-1}}{s}G \rho^{s+1} ,
\formula{frm:stima-generatrice}
$$
with $b\geq0$ and $G>0$.
Then there exists a $\rho^*$ such that if $\rho\leq\rho^*$ then the
operator $T_Y$ and its inverse $T_{Y}^{-1}$
define an analytic canonical transformation on the domain
$\Dscr_{\rho/2}$ .
Moreover we have
$$
\biggl|T_Y z - \sum_{s=0}^{r} E_s z \biggr|_\rho < C^r \rho^{r+2}\ ,
\quad
C=C_1 \max(rB, G)\ ,
$$
with some positive constant $C_1$.
\endclaim

This is a general statement concerning the analyticity of the Lie
transform.  The proof requires minor
changes with respect to similar proofs given, e.g.,
in~\dbiref{Giorgilli-1985} and~\dbiref{Giorgilli-2003}.

In our case the generating sequence $Y$ is replaced by
$Z=\{Z_1,\ldots,Z_r,Q_{r+1},\ldots\}$, and the map is precisely the
coordinate transformation composed with the rotation $\Rmatw$. Remark
that $\Rmatw$ has a trivial impact on the norm.  On the other hand,
from~\frmref{sielie.71} and \frmref{sielie.72} we have
$$
\bignorma{Z_s} \leq T_{s-1} \frac{B_r^{s-1} A}{s}
\leq \frac{1}{\alpha_r^{s-1}} \cdot \frac{B_r^{s-1} A}{s}\ ,\quad
s\leq r\ .
$$
With a trivial estimate we also get
$$
\bignorma{Q_s} \leq \frac{1}{\alpha_r^{s-1}} \cdot \frac{B_r^{s-1} A}{s}\ ,\quad
s> r\ .
$$
This because $Q_s$ corresponds to the right member of~\frmref{map.20b}
with however the sum over $j$ running up to $r$ only.  Therefore,
using also the Diophantine condition, we apply the lemma above
replacing $b$ with $r^{\tau+1} b_1$, with some constant $b_1$
independent of the frequencies $\omega$. With a straightforward use of
Stirling formula, this produces the estimate~\frmref{frm:stimaresto}.

\section{sec:numerica}{Numerical application}

The aim of this section is to present an application of the previously
introduced theory, in the framework of the {\corsivo control of
(simplified models of) particles accelerators}. Our choice has been
motivated by the presence of a large literature (see,
e.g., \dbiref{HanCar-1984}, \dbiref{WanCarSha-1996}, \dbiref{WanCar-2001}, \dbiref{BorCarSkoPapVit-2012}
and \dbiref{BorCarSkoVit-2012}) and by the fact that particles
accelerators are the right benchmarks to apply the Hamiltonian
control, both because the dynamics can be, in very good approximation,
described by a conservative system, and secondly because one can
associate the controller determined by the theory with (a combination
of) basic elements, multipoles that in principle could be inserted in
the accelerator and thus increase the dynamical aperture, namely the size of the domain of (effective) stability of the nominal trajectory.

For sake of simplicity we decided to present our theory for a $2D$
symplectic map that can be used to model a flat beam, that is when the
vertical extent of the beam is much smaller than the horizontal one
(see, e.g., \dbiref{Raubenheimer-1993}, \dbiref{BouSko-2006}
and \dbiref{BorCarSkoPapVit-2012}). The extension to standard $4D$
symplectic maps is straightforward but it involves more complicated
computations.

Let $x_1$ be the horizontal coordinate measured from the nominal circular orbit and $x_2$ its associated momentum, then the beam dynamics can be described by the map:
$$
\pmatrix{x_1' \cr x_2'\cr} = 
\pmatrix{\cos\omega_1& -\sin\omega_1\cr
\sin\omega_1& \cos\omega_1\cr}
\pmatrix{x_1\cr x_2 + x_1^2\cr} 
\formula{map:2D}
$$

The latter map is represented in the form $z' = \Rmatw\circ T_{V} z$
given by~\frmref{map.5} with the generating sequence
$$
V_1 = \pmatrix{0\cr x_1^2\cr}\ ,\quad V_s=0 \quad \forall s\geq 2\ .
$$
Actually in our construction we use the generating sequence $W
= \Rmatw V$.  We remark that for a generic map the generating sequence
$W$ may be easily determined by algebraic manipulation.  We illustrate
the results in a complete form for the value $\omega_1
= \pi(\sqrt{5}-1)$.  A similar calculation has been performed also for
$\omega_1 = \sqrt{2}$ (the choice of $\omega_1 = \sqrt{2}$ instead of
the more usual one $\omega_1 = 2\pi\sqrt{2}$ allows to highlight the
impact of low-order near-resonances), and is illustrated in synthetic
form in section~\sbsref{sbs:sqrt2}.

\subsection{sec:numerica}{Numerical evaluation of the dynamical aperture}
We start our analysis by computing the evolution of a set of points in
a neighborhood of the origin in order to identify the region where
the orbits remained confined for a long but finite time (i.e., the
dynamical aperture).  To this end we take a square of suitable size
filled with 16000 initial points uniformly spaced, and calculate the
evolution for $10^5$ iterations. Then we plot the points that remain
inside the square.  The number $10^5$ is suggested by looking at the
results in~\dbiref{BorCarSkoPapVit-2012}.  The results are reported in
figure~\figref{fig:fignew}.  Our aim is to determine the dynamical
aperture using normal form.  This will likely be a useful tool in
order to investigate the improvement of the dynamical aperture thanks
to the control terms.

\figure{fig:fignew}{
\vbox {\openup1\jot\halign{
{#}\hfil
&\hfil{#}
\cr
{\psfig{figure=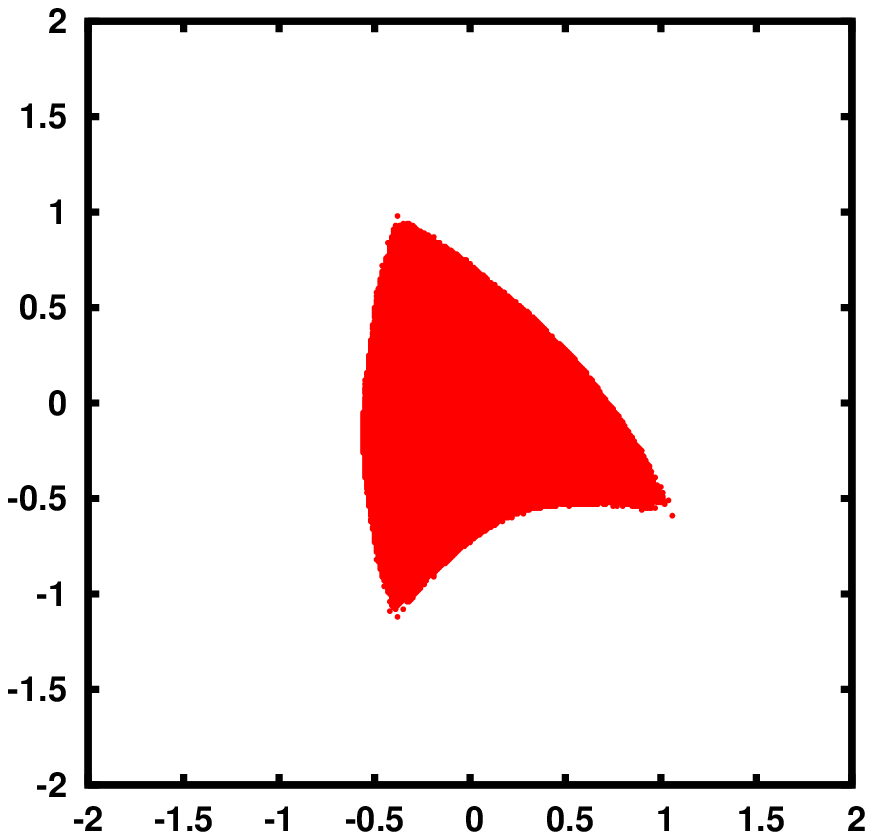,width=\halfwd}}
&{\psfig{figure=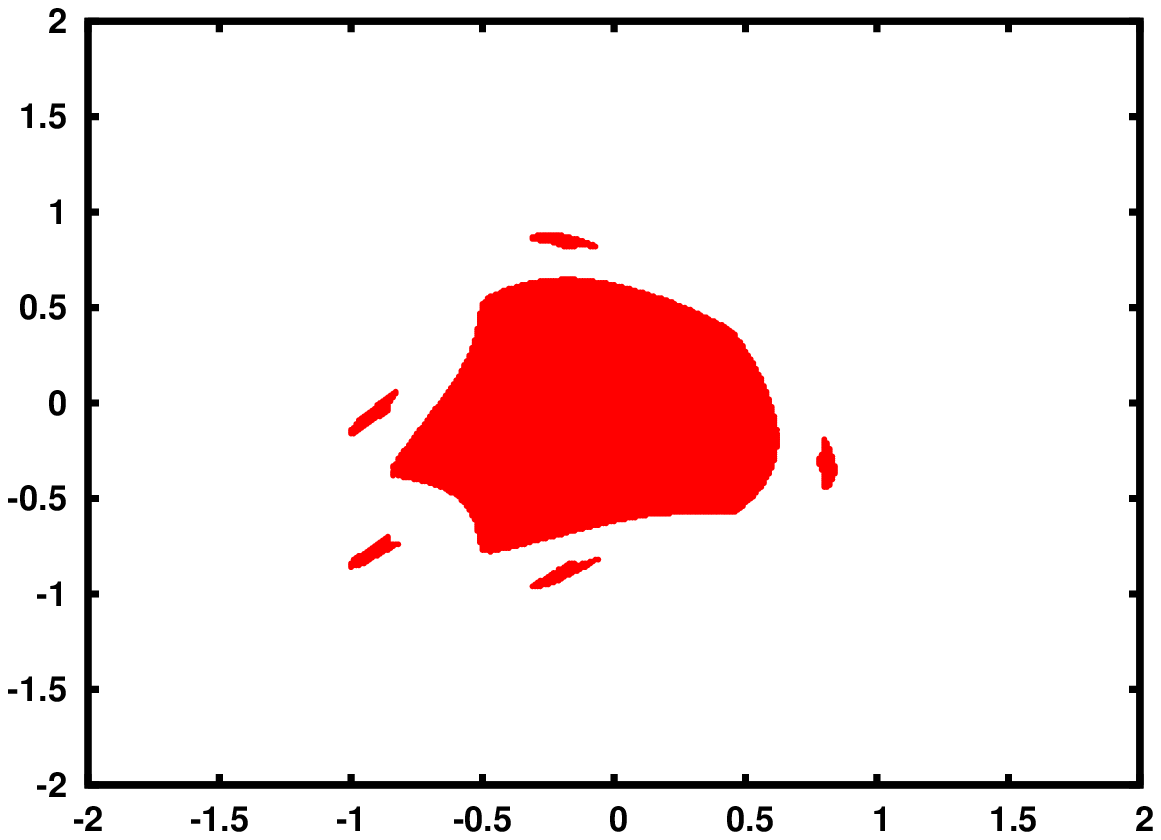,width=\halfwd}}\cr
}}}{Stable region (dynamical aperture) obtained by direct
iteration of the map for $\omega_1 = \pi(\sqrt{5}-1)$ (left-panel) and
$\omega_1 = \sqrt{2}$ (right-panel).}

\subsection{sec:nf}{Construction of the normal form}
The normal form for the map~\frmref{map:2D} is calculated via
algebraic manipulation, representing polynomials with coefficients in
complex form.  Actually we used a package developed on purpose
(see~\dbiref{GioSan-2012}) that we named~\chr.  We perform a number
$r$ of normalization steps.  We note that algebraic manipulation
allows us to take values of $r$ which are relatively large.  E.g.,
$r=50$ or $r=100$ are reachable for a two dimensional map, although
for the present purpose we do stop the normalization at $r=20$.

Therefore, taking a non-resonant frequency $\omega$, we are lead to
consider a map in normal form $z'=T_{Z}\circ {\Rmatw}z$, where
$Z=\{Z_1,\ldots,Z_r,Q_{r+1},\ldots\}$.  Here
$\Qscr=\{0,\ldots,0,Q_{r+1},Q_{r+2},\ldots\}$ is a not normalized
remainder.

Consider now the functions $I_j = (x_j^2+y_j^2)/2$.  They transform as
$$
I' = T_{Z}\circ \Rmatw I = \sum_{s\geq0} E_s^{Z} \Rmatw I\ ,
$$
with the operator $E_s^{Z}$ defined by~\frmref{trslie.2}.  Here we
remark that $\Rmatw I=I$.  On the other hand, in view of
$Z$ being in normal form, we have $E_1^Z I = \ldots = E_r^Z I = 0$,
because $L_{Z_1} I = \ldots = L_{Z_r} I = 0$.  We conclude that for
every integer $t>0$ we have
$$
|I(t) - I(t-1)| \leq \Bigl|\sum_{s>r} E_s^{Z} I(t-1)\Bigr|\ .
\formula{eq:diffusione}
$$
The sum on the r.h.s.~is convergent for every $r$, with a convergence
radius $\rho\to0$ with increasing $r$.

In an analytic approach we could try to estimate a domain of
exponential stability in the spirit of Nekhoroshev's theorem.
However, the analytic estimates, even with the support of algebraic
manipulation, turn out to be too pessimistic and actually unpractical.
Therefore we look for a heuristic estimate of the stability region
in the sense that the orbits are confined in a neighborhood of the
origin for a long time. A similar numerical approach has been previously used in~\dbiref{BorCarSkoPapVit-2012} and~\dbiref{BorCarSkoVit-2012} by computing the fraction of orbits starting in a sphere of a given radius and evolving for a long but finite time without leaving an {a priori} fixed larger sphere.

\subsection{sec:curve}{Level curves of the invariant function}
We enter now the investigation of the invariant function $I(x,y)$
considering different truncation orders $r$ of the normal form.  This
will bring into evidence the asymptotic character of the series, i.e.,
the fact that the radius of the convergence region of the normal form
shrinks to zero for increasing $r$.

\figure{fig:fig2}{
\vbox {\openup1\jot\halign{
{#}\hfil
&\hfil{#}
\cr
{\psfig{figure=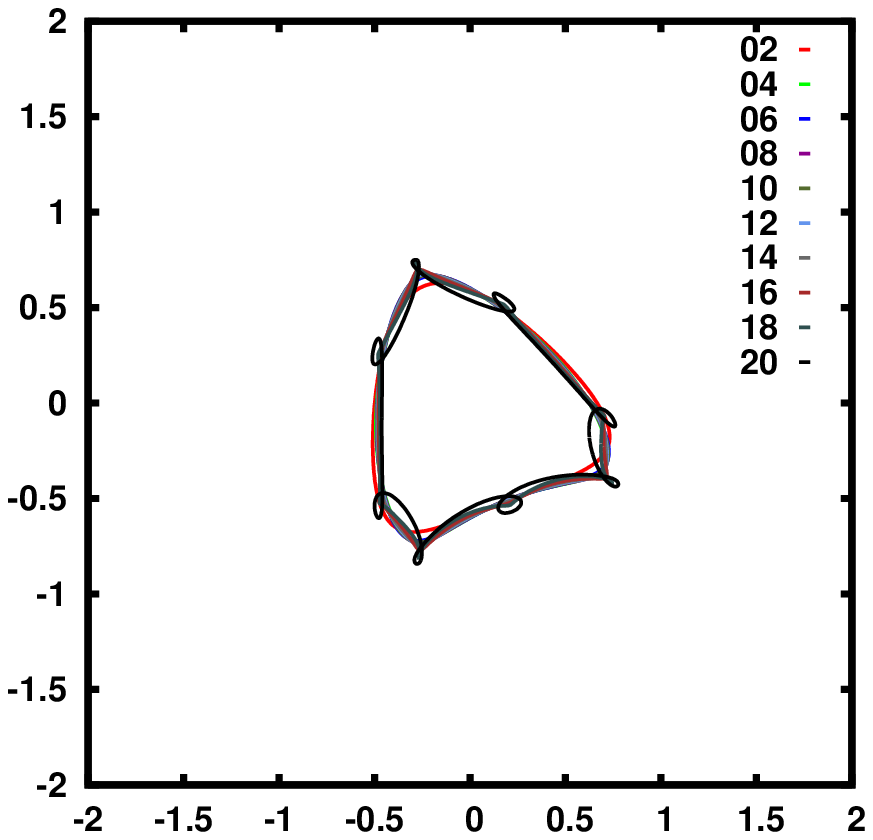,width=\halfwd}}
&{\psfig{figure=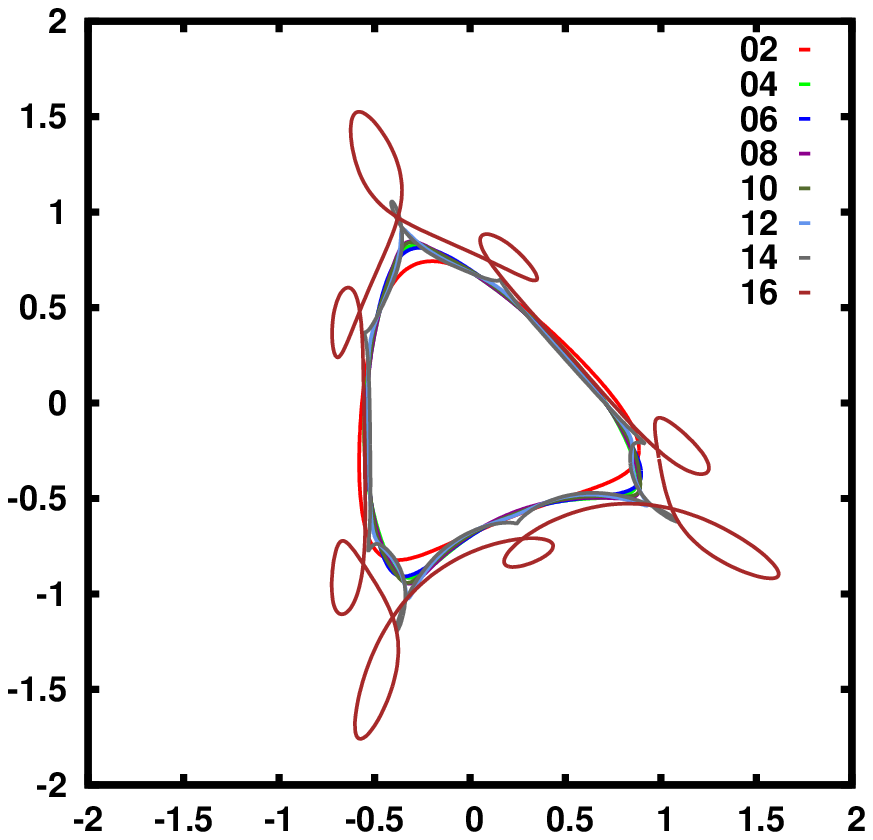,width=\halfwd}}
\cr
{\psfig{figure=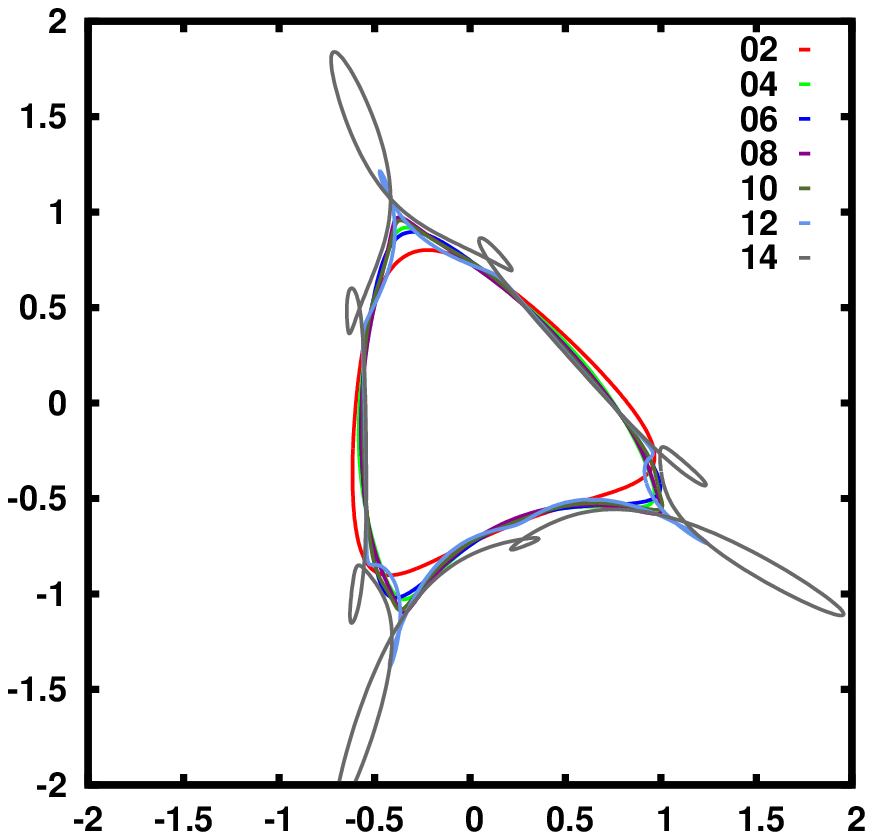,width=\halfwd}}
&{\psfig{figure=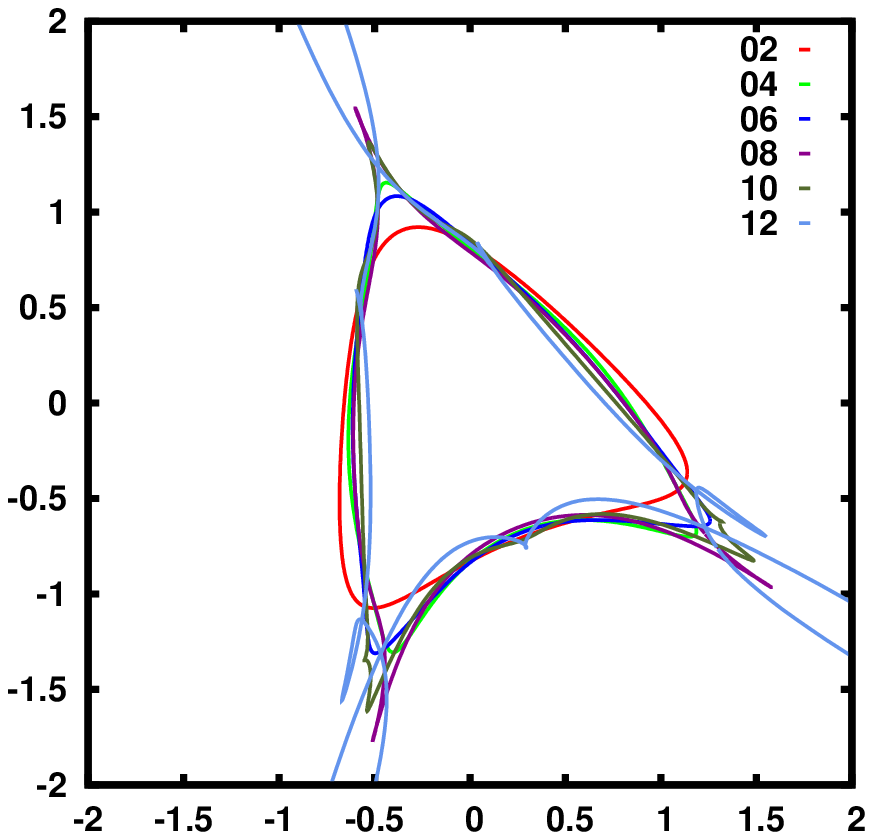,width=\halfwd}}
\cr
}}}{Level curves of the invariant function for the uncontrolled map,
for different values of $\rho=0.60,\, 0.70,\, 0.75,\, 0.85$.  The numbers in the legend correspond to the normalization order $r$ (see text).}

\figure{fig:fig3}{
{\psfig{figure=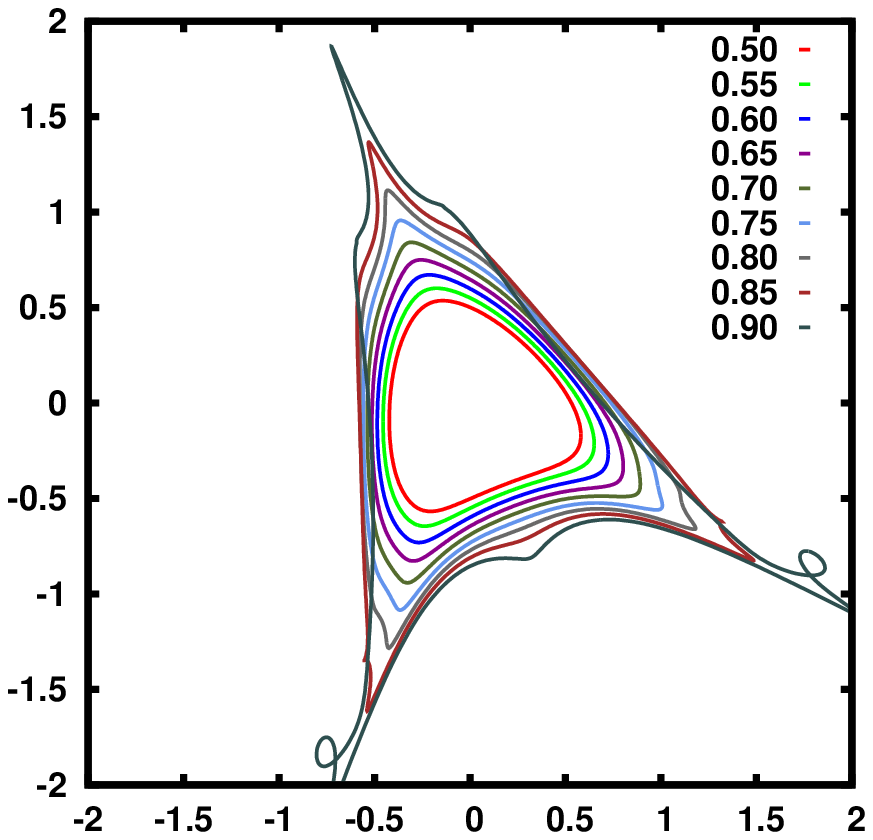,width=0.7\pagewidth}}
}{Level curves of the invariant function for the uncontrolled map for
the normalization order $r=10$ and for different values of $\rho$ as
indicated in the legend.}

\figure{fig:fig4}{
\vbox {\openup1\jot\halign{
\hfil{#}
&{#}\hfil
\cr
{\psfig{figure=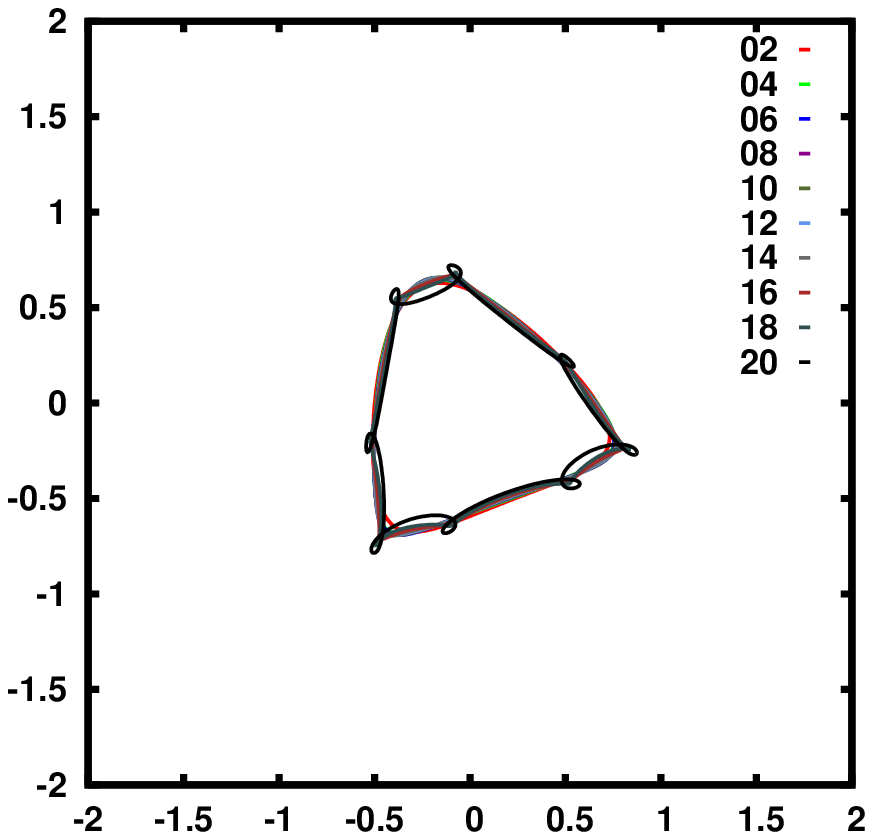,width=\halfwd}}
&{\psfig{figure=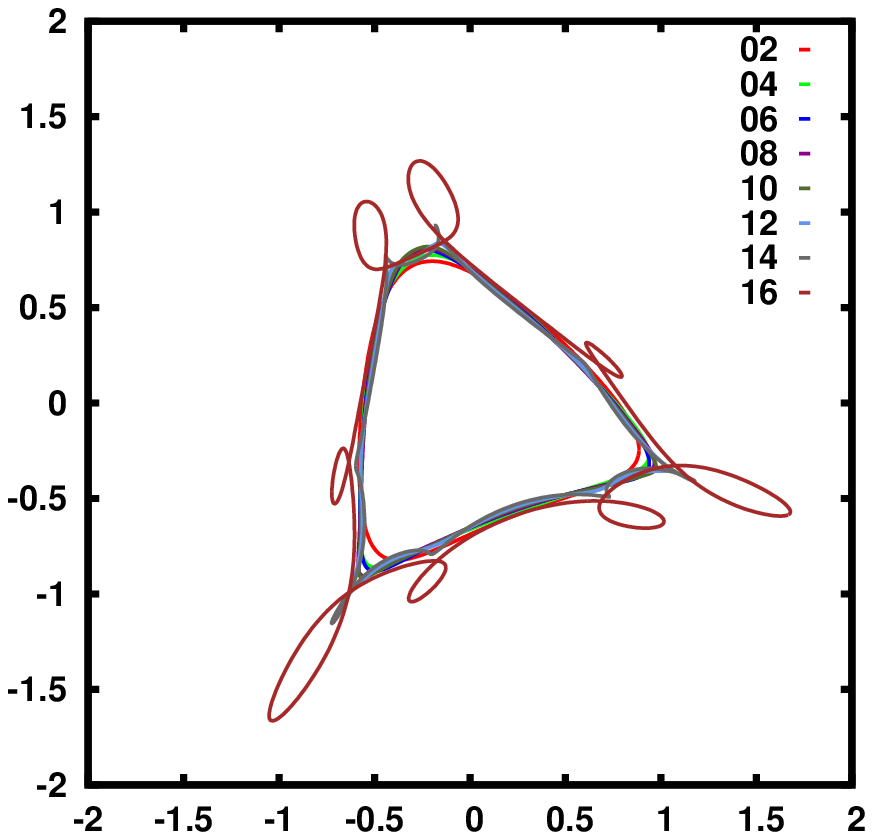,width=\halfwd}}
\cr
{\psfig{figure=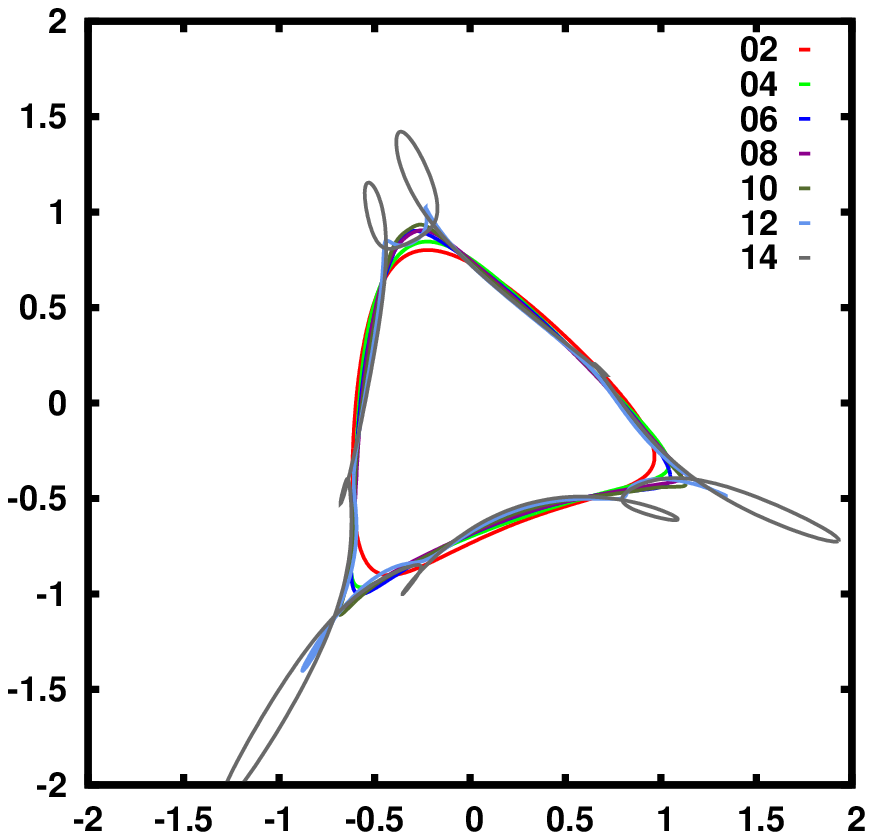,width=\halfwd}}
&{\psfig{figure=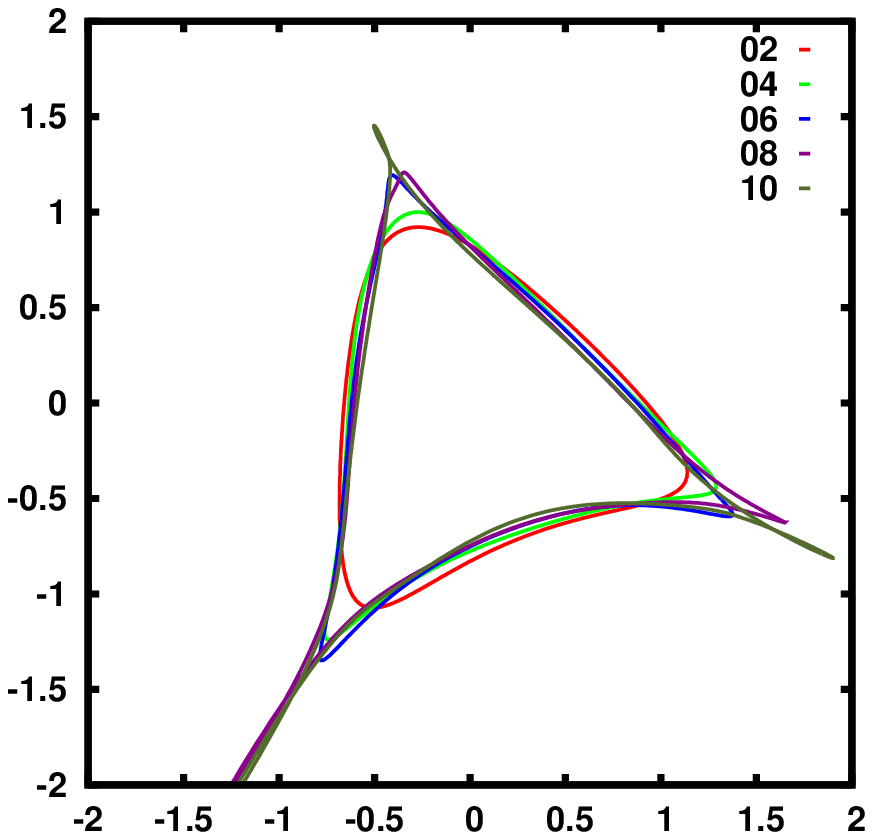,width=\halfwd}}
\cr
}}}{Level curves of the invariant function for the map with the control term $\Fscr_2$,
for different values of $\rho=0.60,\, 0.70,\, 0.75,\, 0.85$.  The numbers in the legend correspond to the normalization order $r$ (see text).}

\figure{fig:fig5}{
{\psfig{figure=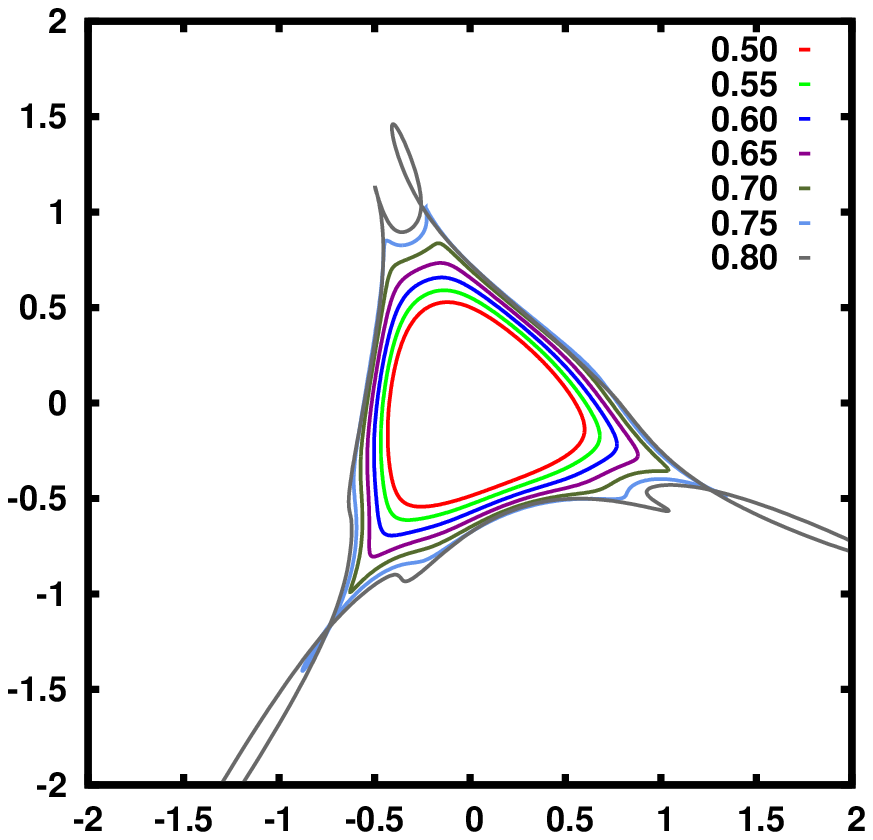,width=0.7\pagewidth}}
}{Level curves of the invariant function for the controlled map for
the normalization order $r=12$ and for different values of $\rho$ as
indicated in the legend.}

\figure{fig:fig6}{
\vbox {\openup1\jot\halign{
{#}\hfil
&\hfil{#}
\cr
{\psfig{figure=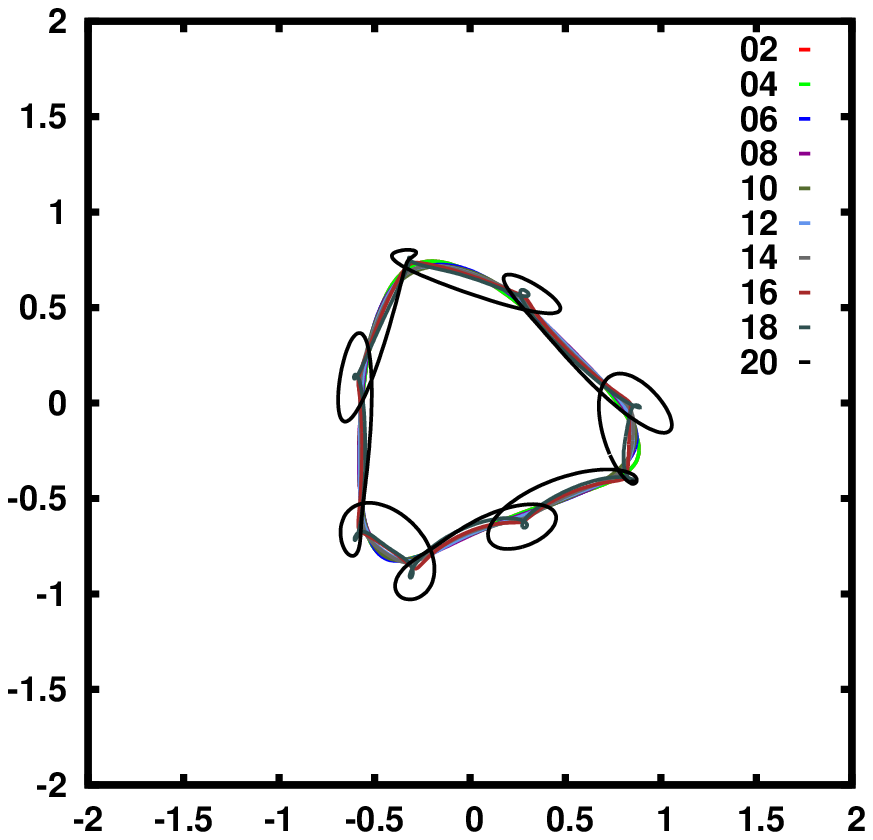,width=\halfwd}}
&{\psfig{figure=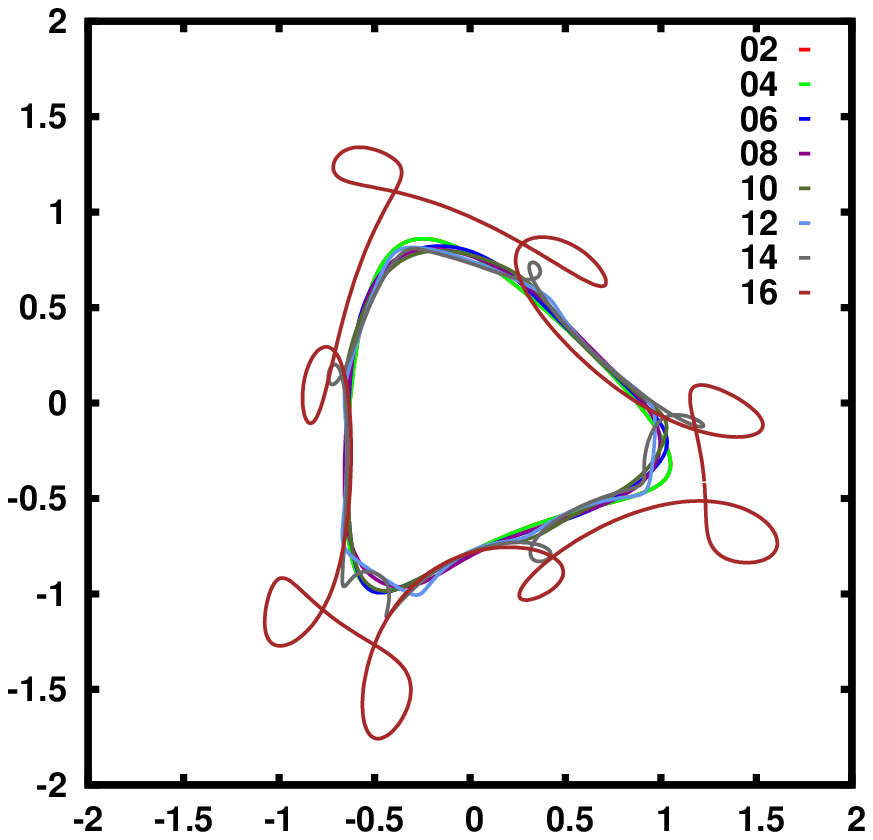,width=\halfwd}}
\cr
{\psfig{figure=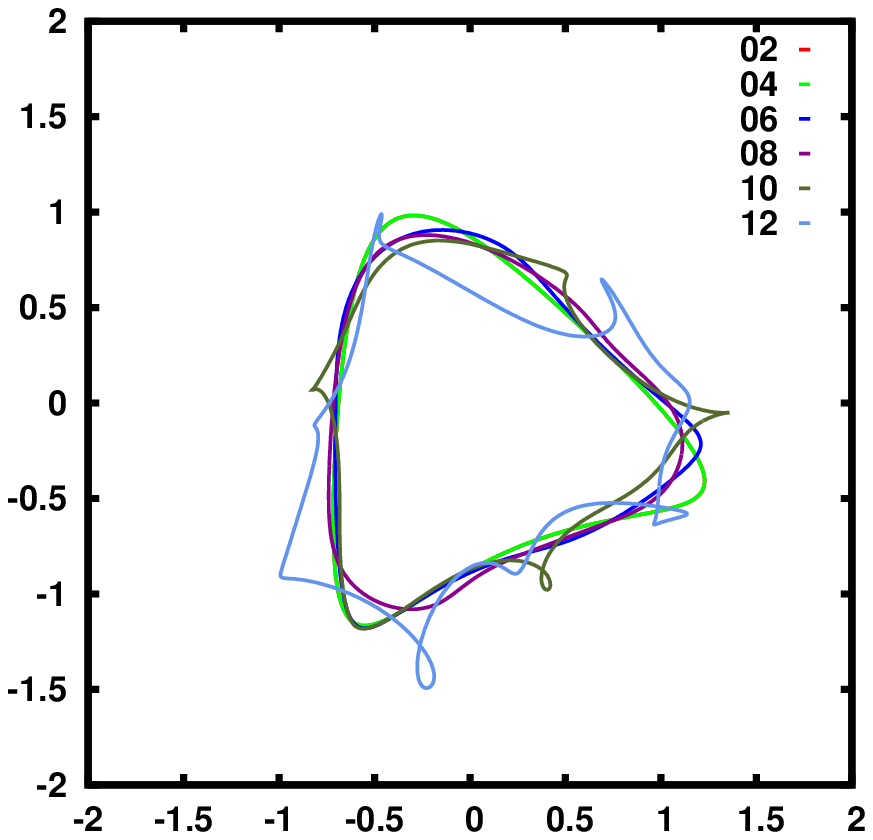,width=\halfwd}}
&{\psfig{figure=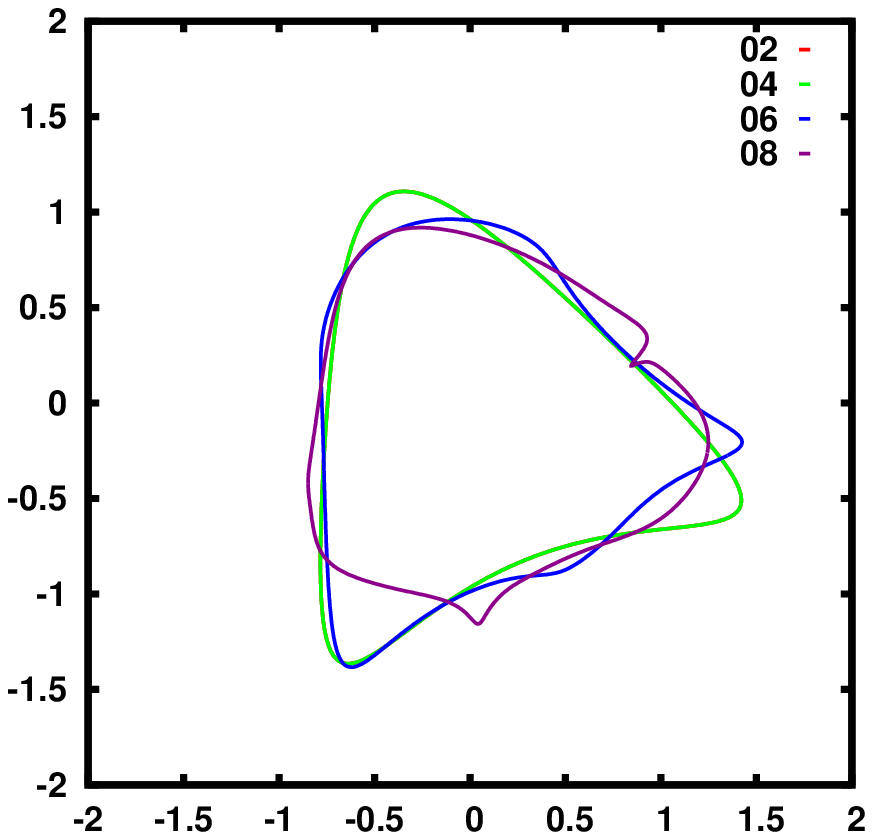,width=\halfwd}}
\cr
}}}{Level curves of the invariant function for the map with the
control term $\Fscr_3$, for different values of $\rho=0.70,\, 0.80,\,
0.90,\, 1.00$.  The numbers in the legend correspond to the
normalization order $r$ (see text).}

\figure{fig:fig7}{
{\psfig{figure=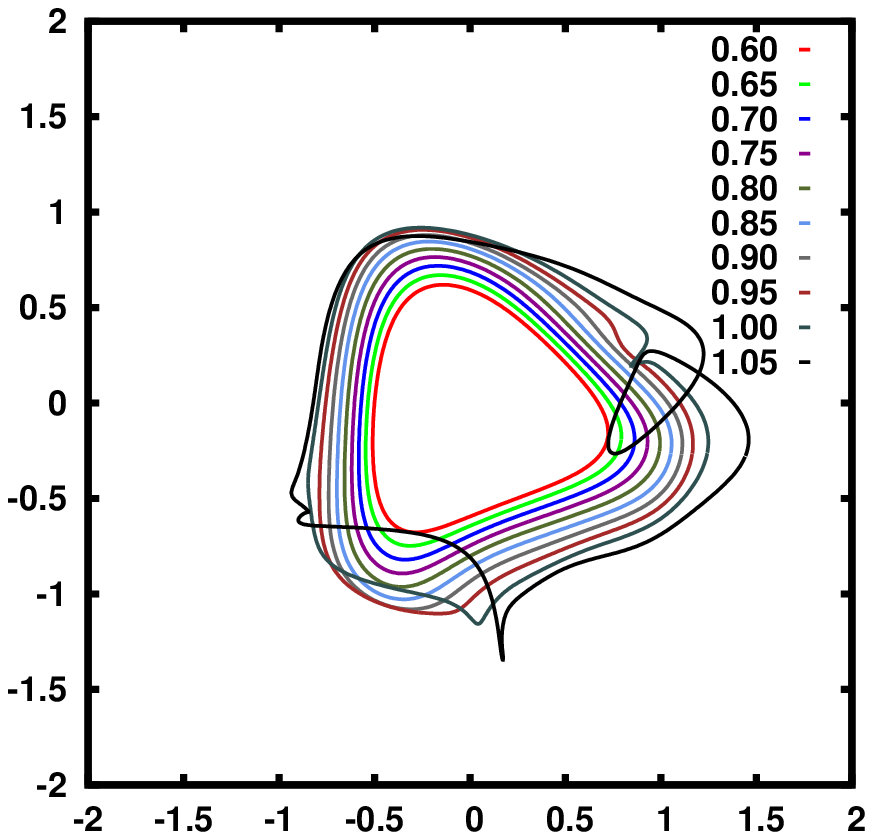,width=0.7\pagewidth}}
}{Level curves of the invariant function for the controlled map for
the normalization order $r=8$ and for different values of $\rho$ as
indicated in the legend.}

\figure{fig:fig8}{
\vbox {\openup1\jot\halign{
{#}\hfil
&\hfil{#}
\cr
{\psfig{figure=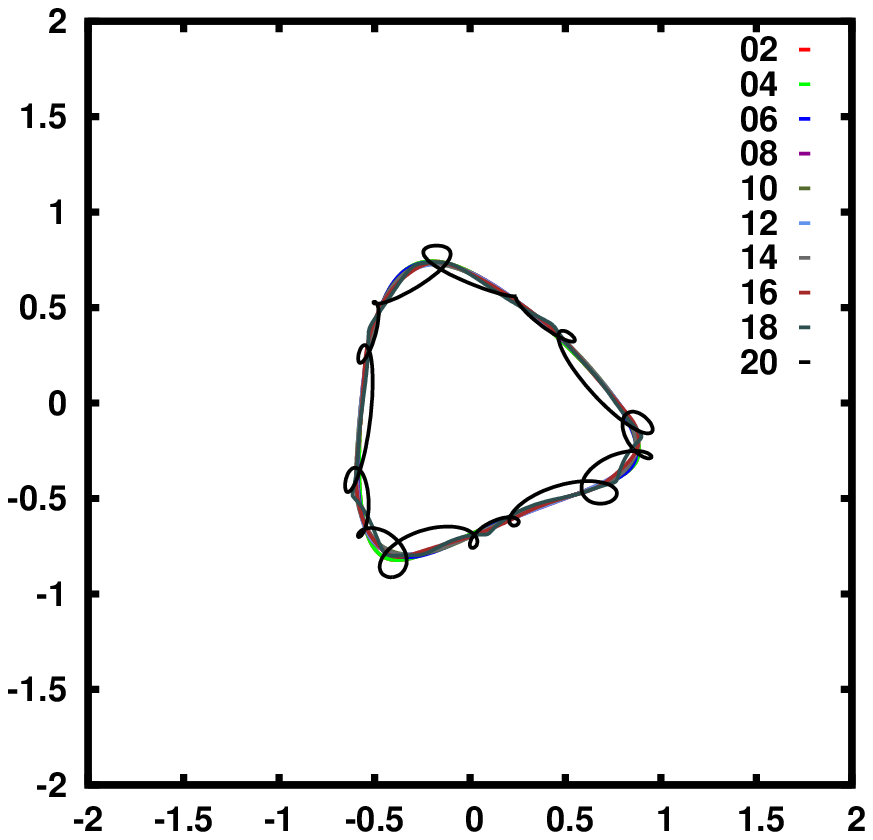,width=\halfwd}}
&{\psfig{figure=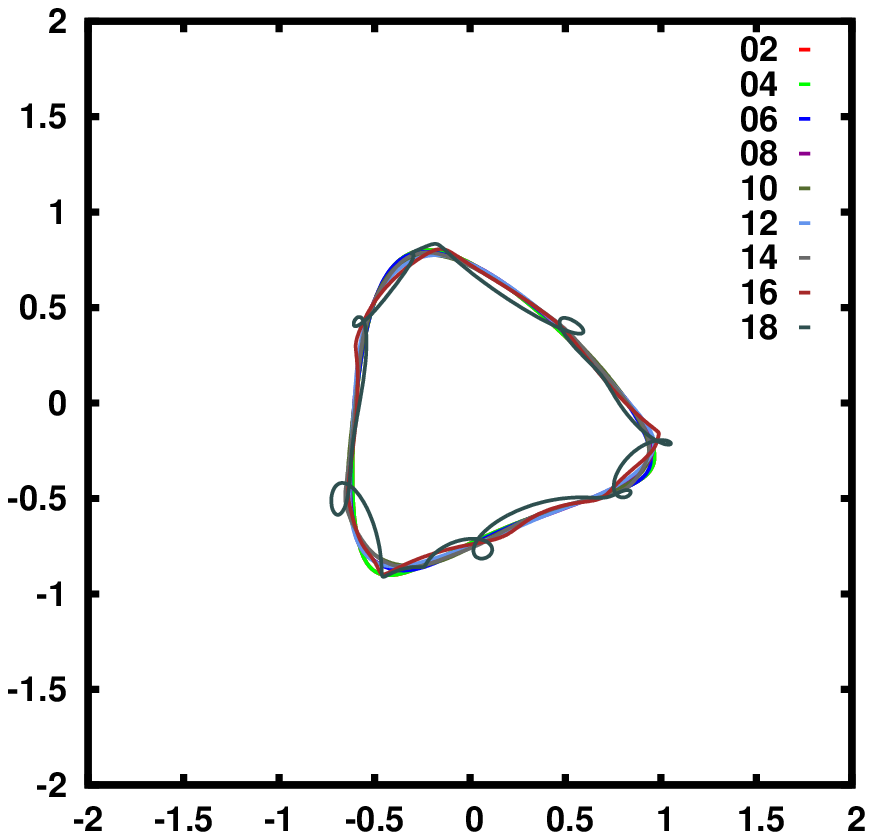,width=\halfwd}}
\cr
{\psfig{figure=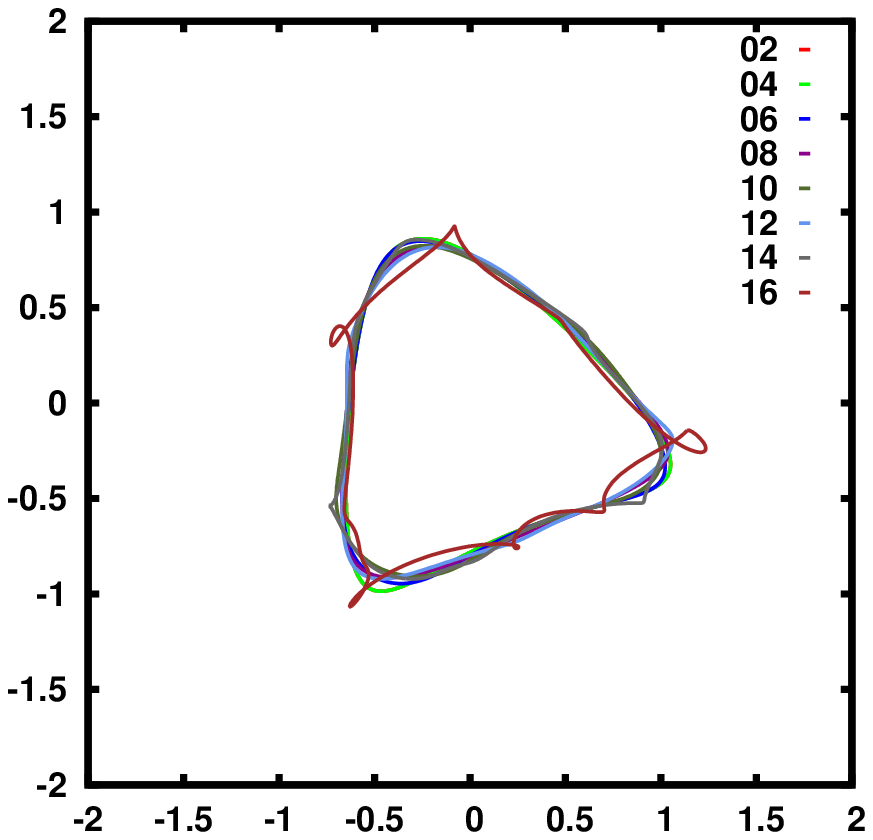,width=\halfwd}}
&{\psfig{figure=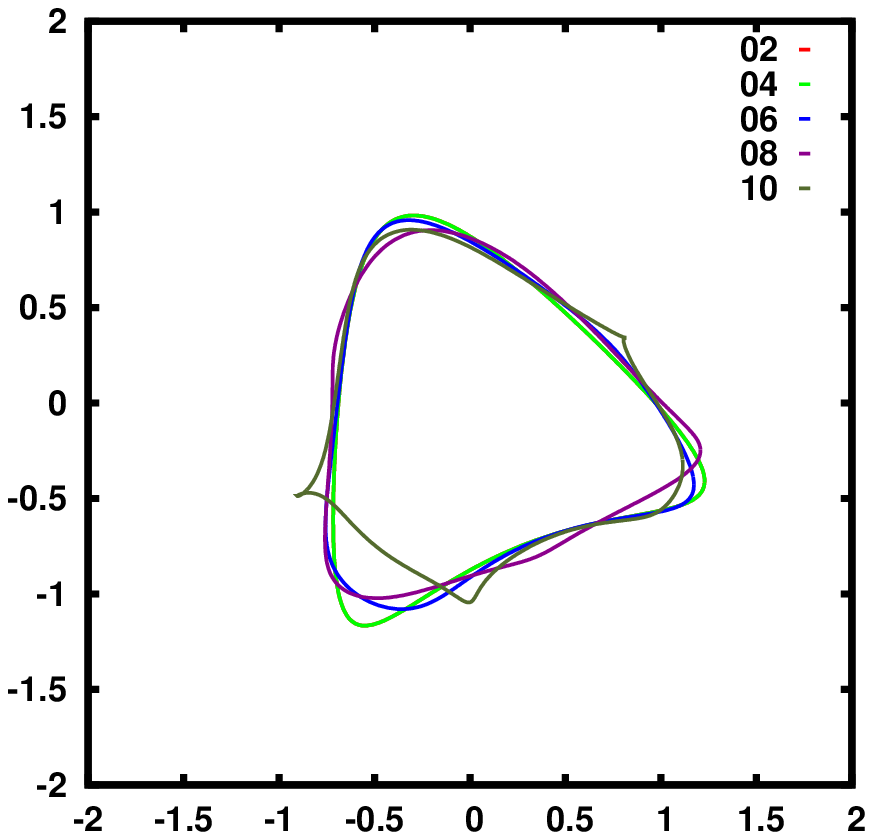,width=\halfwd}}
\cr
}}}{Level curves of the invariant function for the the map with the control term $\Fscr_4$,
for different values of $\rho=0.70,\, 0.75,\, 0.80,\, 0.90$.  The numbers in the legend correspond to the normalization order $r$ (see text).}

\figure{fig:fig9}{
{\psfig{figure=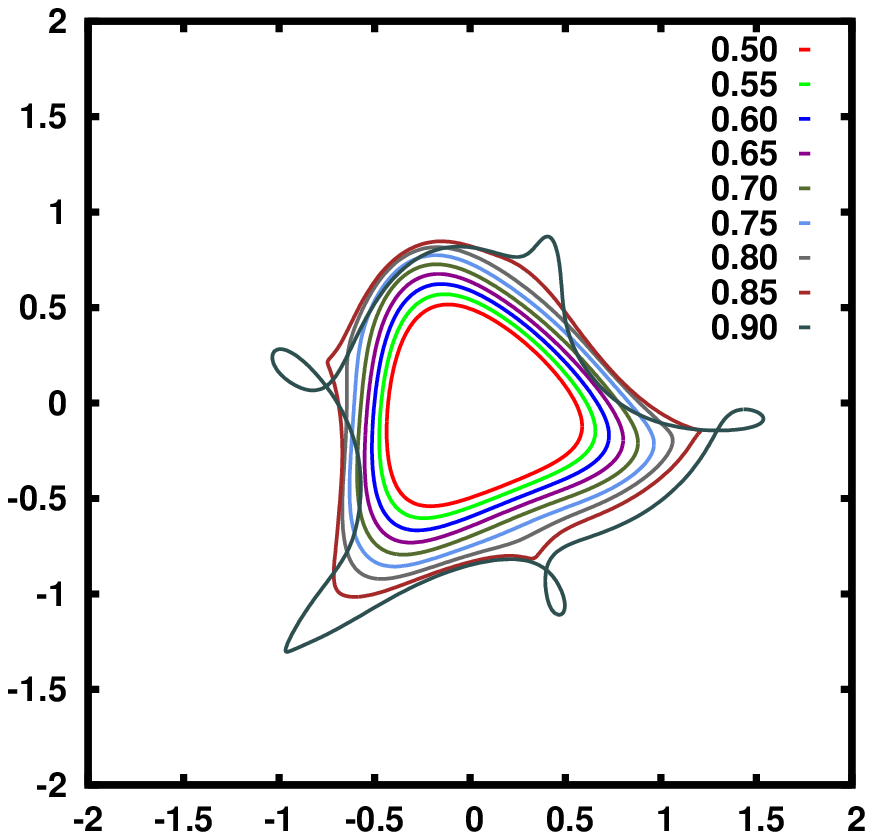,width=0.7\pagewidth}}
}{Level curves of the invariant function for the controlled map for
the normalization order $r=12$ and for different values of $\rho$ as
indicated in the legend.}

We remark that in the normal form coordinates $(x', y')$, the
invariant function is $I=({x'}^2 + {y'}^2)/2$, therefore the invariant
curves are just circles.  In order to compare the curves with the
figures given by the numerical evolution, we apply the transformation
$z=T_{X^{(r)}} z'$, the upper label $r$ denoting the order of
normalization.  According to the asymptotic character of the series,
we expect that for a given radius $\rho$ the image of the circle in
the original coordinates will be a regular curve, just a deformation
of a circle, up to a certain order $r(\rho)$ depending on $\rho$,
while at higher order the image will exhibit a more or less strange
behavior.  Since we have two free parameters, $r$ and $\rho$, we draw
separate figures for fixed $\rho$ and for fixed $r$.  In each panel of 
figure~\figref{fig:fig2} we plot the level curves for a fixed value of
$\rho$ and increasing normalization order $r$.  The four panels
correspond to the values $\rho=0.60,\, 0.70,\, 0.75,\, 0.85$ (see
caption).  The values of $\rho$ are chosen so as to bring into
evidence the asymptotic character of the series.  In the upper-left
panel, with the lowest value of $\rho$, we see that most curves
corresponding to orders $r=2,\,\ldots,\,18$ seem to visually coincide,
as expected from a convergent series.  However, at $r=20$ we see that
the curve exhibits a strange behavior, forming some unexpected loops.
The curves would become more and more tangled at higher orders, data
not shown.  We may say that there is an interval of apparent
convergence, in this case up to $18$.  A similar behavior shows up
for $\rho = 0.70$ (upper-right panel) and $\rho= 0.75$ (lower-left
panel), but the interval of apparent convergence is more and more
restricted.  E.g., for $\rho=0.70$ the curves corresponding to
$r=4,\,6,\,8,\,10$ are quite close and still regular, while a loop
appears at $r=12$.  For the higher value $\rho=0.85$ all curves are
well separated, so that the phenomenon of apparent convergence
completely disappears.

In figure~\figref{fig:fig3} we fix $r=10$ and draw the level curves
for increasing values of $\rho$.  The choice $r=10$ is suggested by
the preliminary analysis performed on figure~\figref{fig:fig2} (see
lower-left panel of such figure).  As we see, for $\rho=0.9$ the curve
exhibits a loop, thus showing that we are out of the radius of apparent
convergence.  Thus we may consider the interval $I_{\rho}=(0.7{\rm -}0.8)$
as the right one, and we may make a safe choice with ${\rho}=0.75$.

The comparison between the two figures provides us with a heuristic
criterion in order to make a safe choice of $\rho$, which is the
relevant quantity for the dynamics.  A corresponding value of $r$ is
also suggested. The aim of the next section is to compare the improved
behavior and the optimality of the controlled map against the
uncontrolled one using such heuristic metrics: ${\rho}$ and $r$.

\figure{fig:fig10}{
\vbox {\openup1\jot\halign{
{#}\hfil
&\hfil{#}
\cr
{\psfig{figure=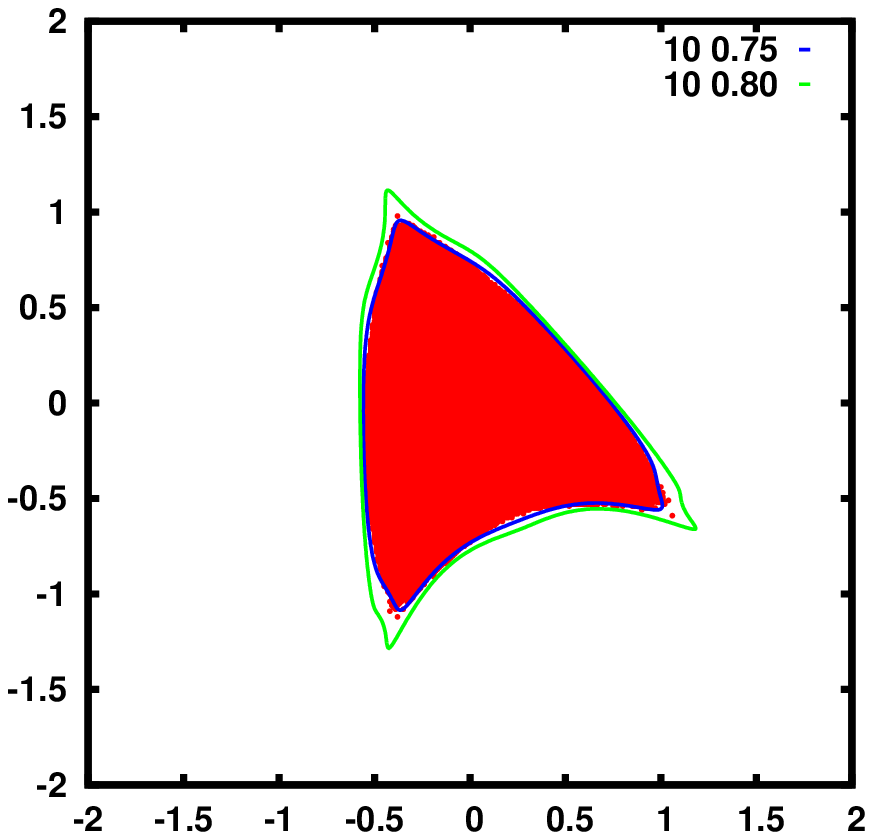,width=\halfwd}}
&{\psfig{figure=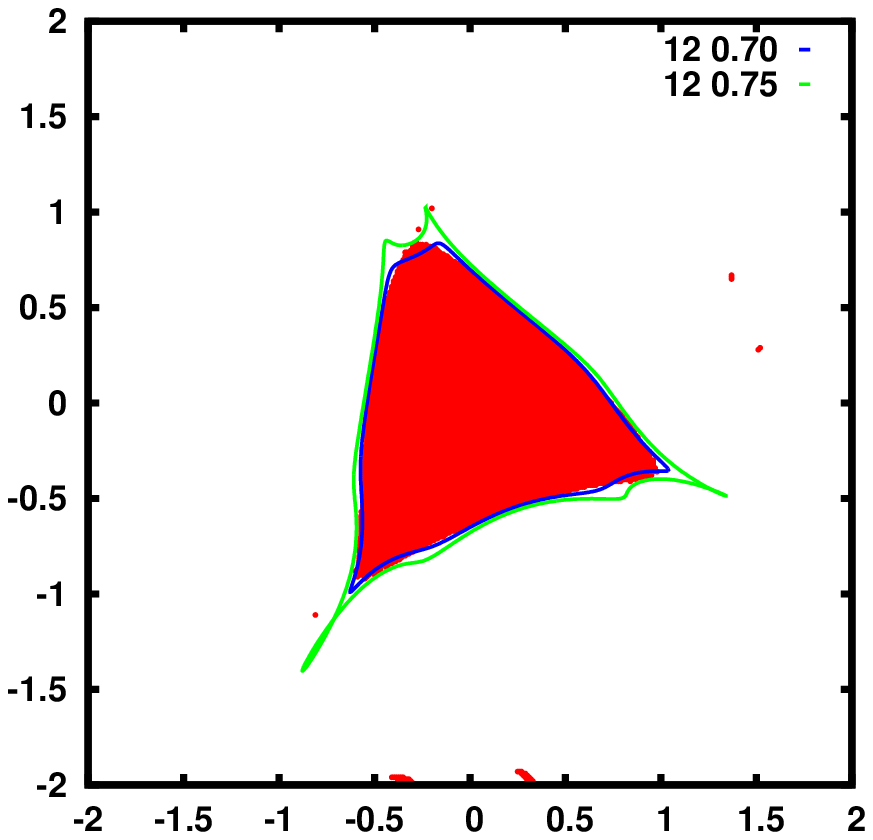,width=\halfwd}}
\cr
 {\psfig{figure=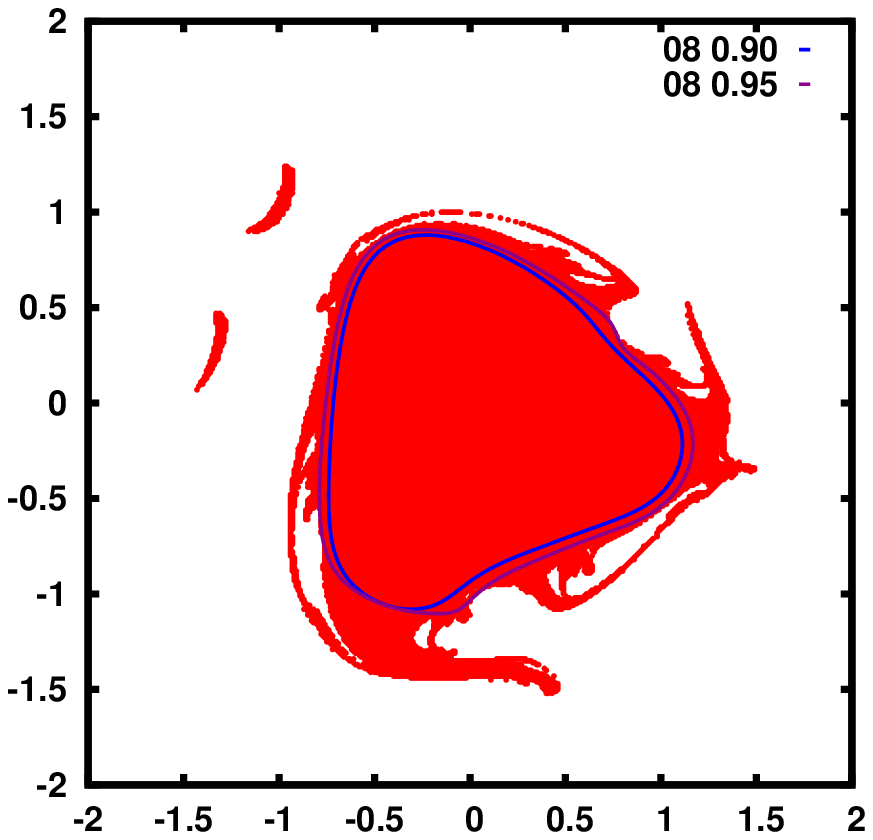,width=\halfwd}}
&{\psfig{figure=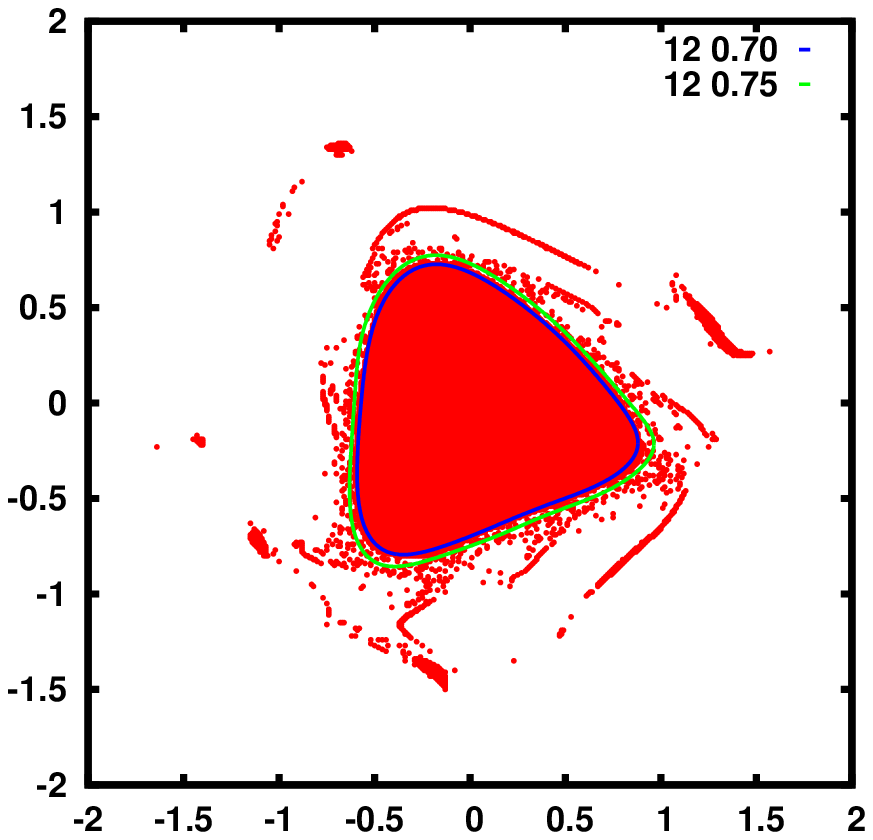,width=\halfwd}}
\cr
}}}{Comparison of the level curves of the invariant functions with the
stability zone obtained by direct iteration of the map (with $\omega_1 = \pi(\sqrt{5}-1)$).  The four
panels (left to right) correspond to the uncontrolled map and the
three different controls $\Fscr_2$, $\Fscr_3$ and $\Fscr_4$.  The
normalization order $r$ and the radius $\rho$ are selected using the
heuristic criterion illustrated in subsection~\sbsref{sec:curve}.
We plot the level curves for two different values of $\rho$.  The
lower one produces a reliable approximation of the stability zone.
One sees that the best control terms is $\Fscr_3$ (left-lower panel).}

\figure{fig:fig11}{
\vbox {\openup1\jot\halign{
{#}\hfil
&\hfil{#}
\cr
{\psfig{figure=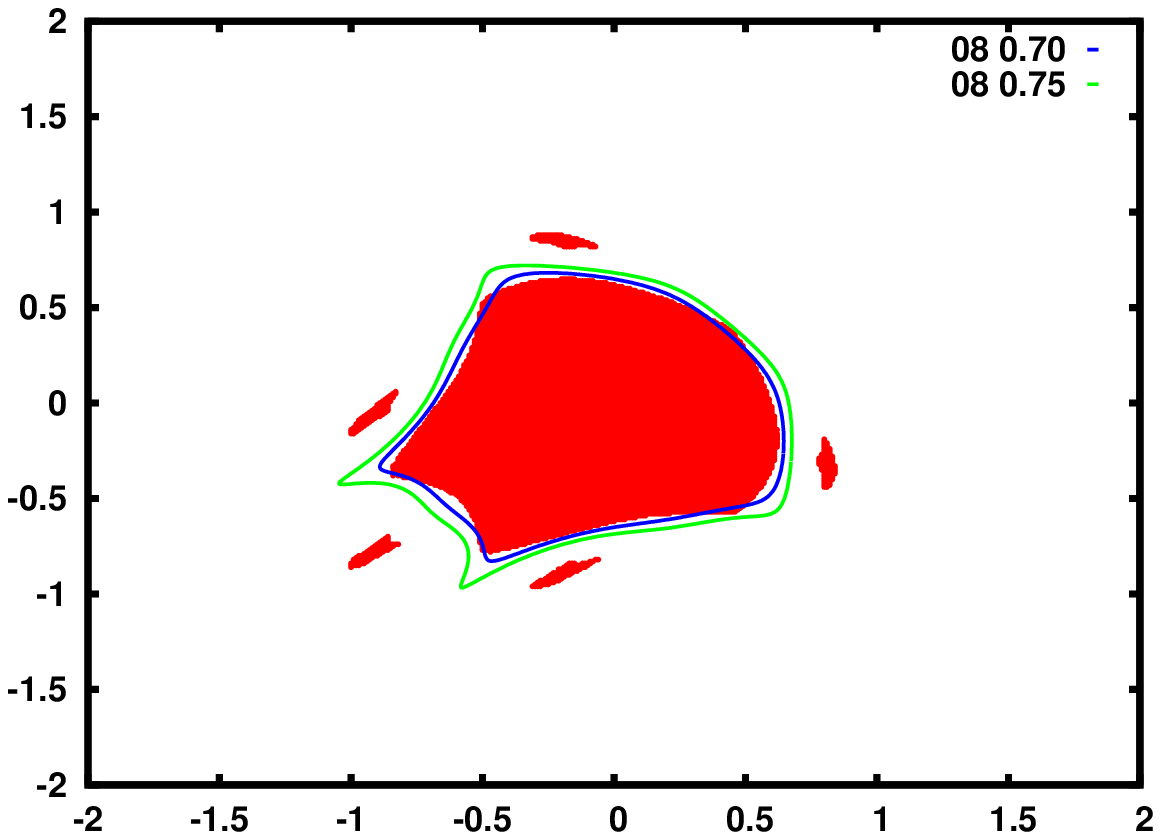,width=\halfwd}}
&{\psfig{figure=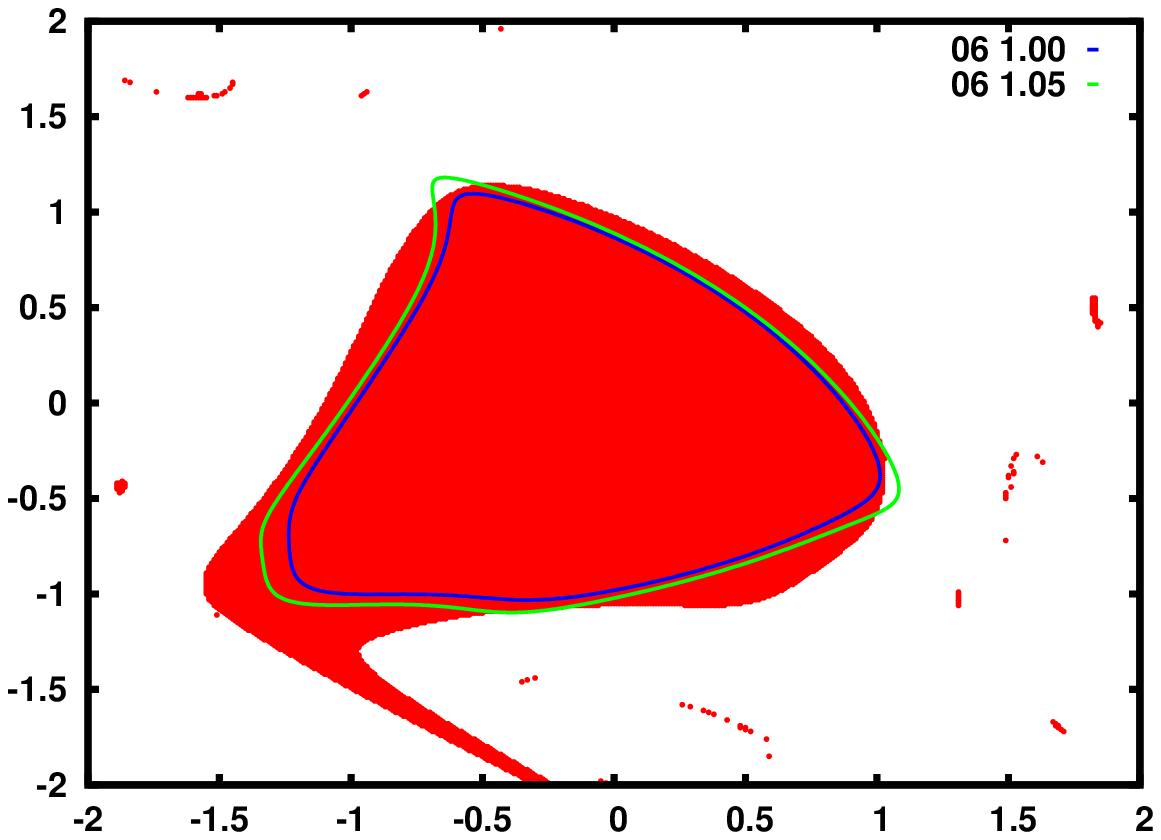,width=\halfwd}}
\cr
{\psfig{figure=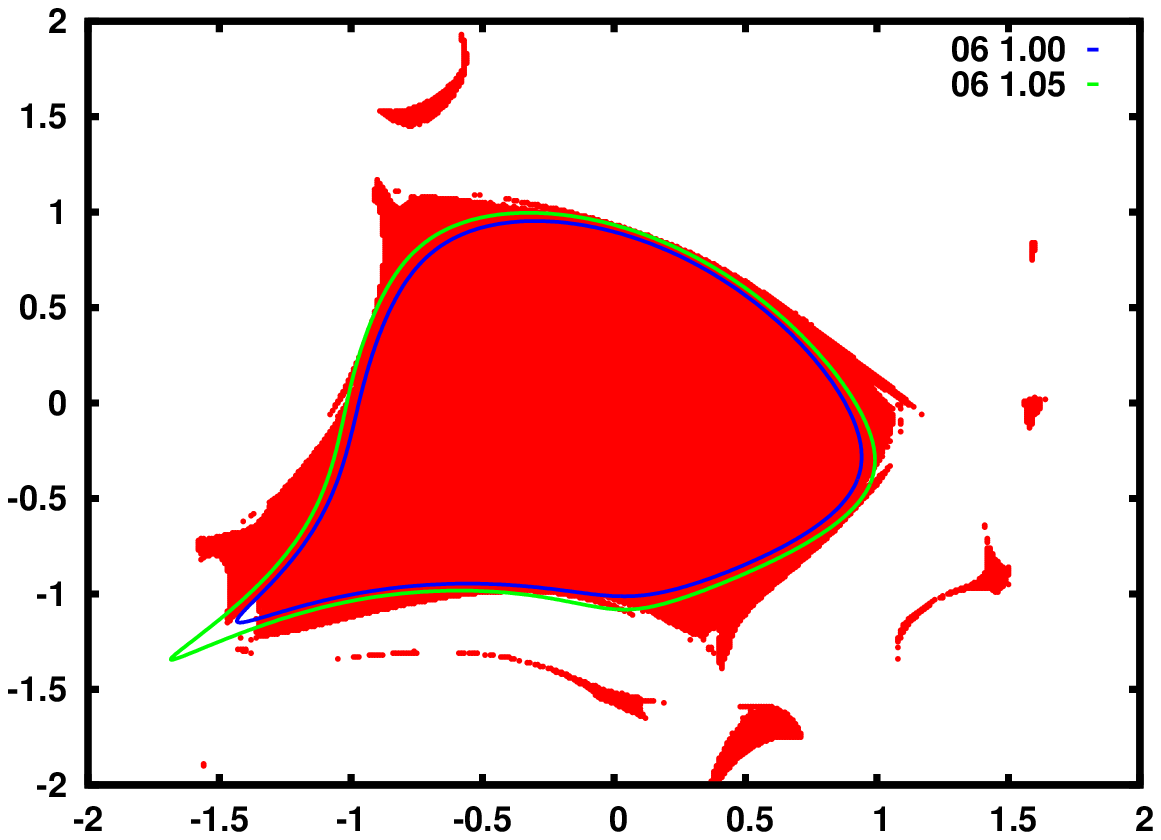,width=\halfwd}}
&{\psfig{figure=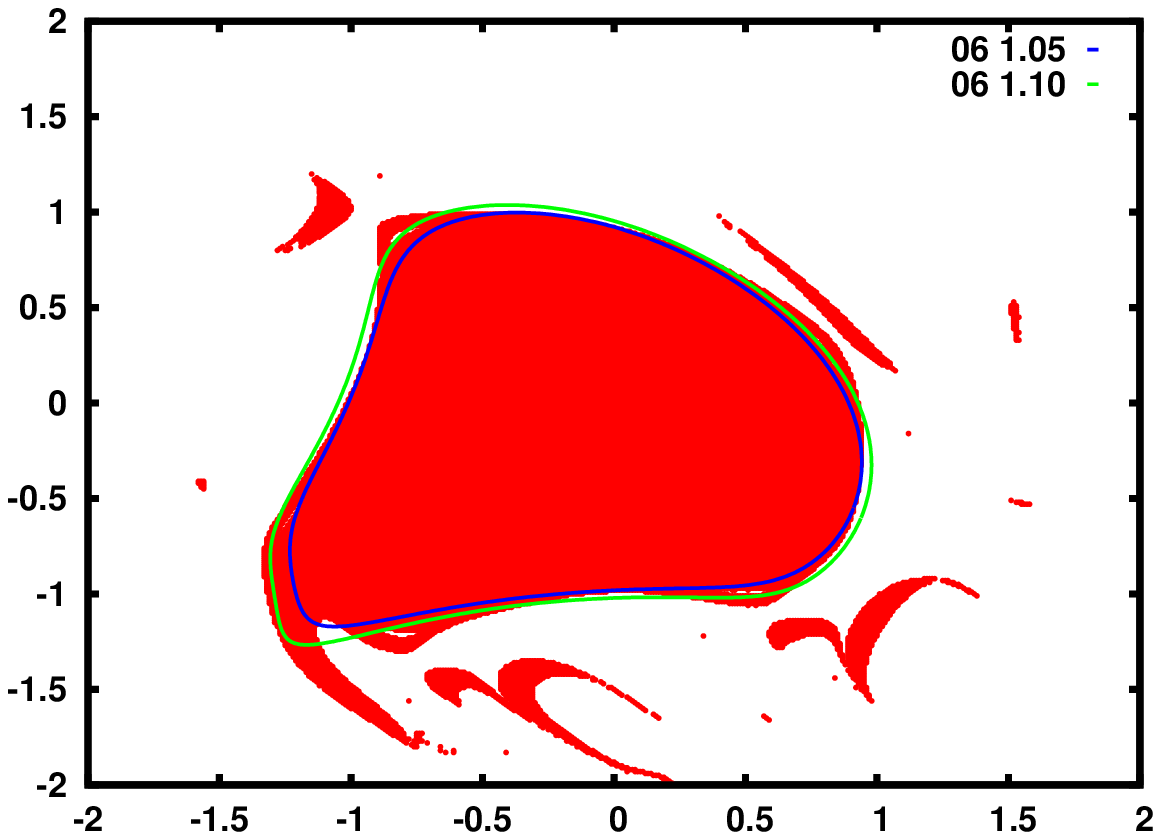,width=\halfwd}}
\cr
}}}{Comparison of the level curves of the invariant functions with the
stability zone obtained by direct iteration of the map (with $\omega_1 = \sqrt{2}$). See caption of figure~\figref{fig:fig10}.}

\subsection{sbs:controllo}{Level curves for the controlled map}
We illustrate the application of the criterion above by introducing
different corrections to the map.

We introduce the controlled map as illustrated in
section~\sbsref{sbs:nfmappe}.  We consider three different cases by
choosing the control terms $\Fscr_2=\{0,F_{2},0,\ldots\}$,
$\Fscr_3=\{0,F_{2},F_{3},0,\ldots\}$ and
$\Fscr_4=\{0,F_{2},F_{3},F_{4},0,\ldots\}$, with $F_{s}$ given
by~\frmref{eq:controllo}.  This corresponds to introducing controls including increasing orders.

The results for the map with the control term $\Fscr_2$ are
illustrated in figures~\figref{fig:fig4} for fixed values of $\rho$
(analogous to figure~\figref{fig:fig2} for the uncontrolled map)
and~\figref{fig:fig5} for a fixed value of $r=12$ (analogous to
figure~\figref{fig:fig3}).  Looking at the latter figure, we may
consider as acceptable a value $\rho$ in the interval $0.65$--$0.75$.

Similar figures are reported also for the control term $\Fscr_3$
in~\figref{fig:fig5}--\figref{fig:fig7} and $\Fscr_4$
in~\figref{fig:fig8}--\figref{fig:fig9}.  In the former case our
heuristic criterion suggests the value $\rho\simeq0.9$, in the latter
case we get $\rho\simeq0.7$.  Thus it seems reasonable to conclude
that the optimal control term, to achieve a good compromise between to have a large $\rho$ but relatively not so large $r$, is $\Fscr_3$.

\subsection{sec:confronto}{Comparison with the dynamics}
We now compare the level curves of the invariant function with the
results of the dynamics obtained by direct numerical iteration of the
map (see subsection~\sbsref{sec:numerica}).

For the uncontrolled map and the three controlled map considered
above, we superimpose to the effectively stable points defined above, the level curves
corresponding to $\rho=0.75$ for the uncontrolled map, $\rho=0.70$ for
the control term $\Fscr_2$, $\rho=0.90$ for the control term $\Fscr_3$
and $\rho=0.70$ for the control term $\Fscr_4$.  The results are
reported in figure~\figref{fig:fig10}.  As one sees there is a good
agreement between the dynamical aperture suggested by the dynamics and
the one estimated via the normal form.  For comparison we also add the
level curves for values of $\rho$ slightly bigger than the good
estimated ones, thus showing that actually our heuristic criterion
produces reliable hints.

\subsection{sbs:sqrt2}{Changing the rotation number}
As anticipated, we performed a complete calculation for the case
$\omega_1=\sqrt{2}$.  The resulting figures for the level curves at
different values of $\rho$ and $r$ exhibit a behavior similar to
figures~\figref{fig:fig2}--\figref{fig:fig10}.  Thus we do not include
a complete set of figures, reporting only the final results in
figure~\figref{fig:fig11} (corresponding to figure~\figref{fig:fig10} of the previous section).

\acknowledgements{A.~G. and M.~S. have been partially supported by the
research program ``Teorie geometriche e analitiche dei sistemi
Hamiltoniani in dimensioni finite e infinite'', PRIN 2010JJ4KPA\_009,
financed by MIUR.  The work of T.C. presents research results of the
Belgian Network DYSCO (Dynamical Systems, Control, and Optimization),
funded by the Interuniversity Attraction Poles Programme, initiated by
the Belgian State, Science Policy Office.}

\references

\bye